\documentclass[a4paper,10pt]{article}

\usepackage{graphicx}%
\usepackage{multirow}%
\usepackage{amsmath,amssymb,amsfonts}%
\usepackage{amsthm}%
\usepackage{mathrsfs}%
\usepackage[title]{appendix}%
\usepackage{xcolor}%
\usepackage{textcomp}%
\usepackage{manyfoot}%
\usepackage{booktabs}%
\usepackage{algorithm}%
\usepackage{algorithmicx}%
\usepackage{algpseudocode}%
\usepackage{listings}%
\usepackage{mathtools}
\usepackage{float}
\usepackage[section]{placeins}
\usepackage{caption}
\usepackage{subcaption}
\usepackage{setspace}
\usepackage{algpseudocode}
\algdef{SE}[IF]{NoThenIf}{EndIf}[1]{\algorithmicif\ #1}{\algorithmicend\ \algorithmicif}

\newcommand{\rrk}{\mathbf{r_k}}
\newcommand{\trk}{\mathbf{\tilde{r}_k}}
\newcommand{\sk}{\mathbf{s_k}}
\newcommand{\tsk}{\mathbf{\tilde{s}_k}}
\newcommand{\rrkm}{\mathbf{r_{k-1}}}
\newcommand{\trkm}{\mathbf{\tilde{r}_{k-1}}}

\newcommand{\twkp}{\mathbf{\tilde{w}'_k}}

\newcommand{\twk}{\mathbf{\tilde{w}_k}}

\newcommand{\wk}{\mathbf{w_k}}

\newcommand{\M}{\mathbf{M^{-1}}}
\newcommand{\A}{\mathbf{A}}

\usepackage[sort&compress,numbers,square]{natbib}
\begin{document}

\title{The Detection and Correction of Silent Errors in Pipelined Krylov Subspace Methods}

\author{Erin Claire Carson\footnote{Charles University, Faculty of Mathematics and Physics, Charles University, Prague, Czech Republic. carson@karlin.mff.cuni.cz, hercik.j@outlook.com} \and Jakub Herc\'{i}k\footnotemark[1]}
\date{}

\maketitle

\begin{abstract}As computational machines become larger and more complex, the probability of hardware failure rises. ``Silent errors'', or bit flips, may not be immediately apparent but can cause detrimental effects to algorithm behavior. In this work, we examine an algorithm-based approach to silent error detection in the context of pipelined Krylov subspace methods, in particular, Pipe-PR-CG, for the solution of linear systems. Our approach is based on using finite precision error analysis to bound the differences between quantities which should be equal in exact arithmetic. By monitoring select quantities during the iteration, we can detect when these bounds are violated, which indicates that a silent error has occurred. We use this approach to develop a fault-tolerant variant and also suggest a strategy for dynamically adapting the detection criteria. Our numerical experiments demonstrate the effectiveness of our approach.\end{abstract}

\section{Introduction}\label{sec:intro}

We consider the problem of solving the linear system $Ax=b$, where $A\in \mathbb{R}^{n\times n}$ is symmetric positive definite (SPD). We are particularly interested in the case where $A$ is very large and sparse, in which case the conjugate gradient method is usually the method of choice.

In large-scale settings, the high cost of communication in the classical conjugate gradient method of Hestenes and Stiefel \cite{OGCG}, which we refer to as HS-CG, motivated the search for mathematically equivalent CG variants better suited for implementation on parallel machines. One possibility is to introduce auxiliary vectors and rearrange the procedure in a way that requires only a single synchronization point, and ``pipeline'' the inner product reductions and matrix operations to occur concurrently, so that global synchronization no longer causes a bottleneck. This is the strategy introduced in the ``pipelined CG'' method of Ghysels and Vanroose \cite{GVCG}. However, such communication-hiding variants may amplify  numerical problems which already exist in HS-CG such as delayed convergence, or may worsen the maximal attainable accuracy (the level at which the approximation error $||\mathbf{x - x_k}||$ starts to stagnate) due to the sensitivity of CG to rounding errors \cite{ChenCarson}.

A solution to these numerical problems was developed in \cite{ChenCarson}. Here, the authors employ the newly introduced expressions as ``predicted'' values, which are then used to compute the next few quantities until they are ``recomputed'' using the original formulas. This so-called ``predict and recompute'' variant (Pipe-PR-CG) helps mitigate the deviation of true quantities and their values approximated by recurrences. This idea was based on the previous work of Meurant \cite{MeurantCray} which aimed to stabilize the HS-CG algorithm while also retaining the potential for parallelism.

As computer hardware grows in scale, the probability of failure grows, and thus the topic of error handling and fault-tolerant algorithms is more important than ever. Some errors are simple to detect, since they, e.g., result in a crash of the computation. However, ``silent errors'' may not be immediately apparent, but can result in significantly altered algorithm behavior. There has been a significant amount of work in studying the effect and detection of silent errors in numerical linear algebra computations; see, e.g., \cite{agullo2017hard, MeurantCG, ElliottIEEE, Aupy2017, elliott2015numerical, bronevetsky2008soft}. One approach is to perform the computation multiple times and compare the results. However, this substantially increases the overall computational cost, in terms of both time and energy. A more efficient approach is desirable.

In this work, we derive an approach for silent error detection in Pipe-PR-CG and develop a modified version of the algorithm which is able to automatically detect and correct the silent faults. Our detection methods are based on the comparison of ``gaps'' between certain quantities which are equal in exact arithmetic and the bounds on their values in finite precision. After providing relevant background in Section \ref{sec:background}, in Section \ref{sec:effect} we provide a set of experiments which demonstrates the sensitivity of Pipe-PR-CG to silent errors. Section \ref{sec:detection} focuses on constructing several criteria for the detection of silent errors in the Pipe-PR-CG algorithm based on floating point rounding error analysis, and in Section \ref{sec:ftpipeprcg}, we present a fault-tolerant variant and experimental results demonstrating the effectiveness of the detection criteria. Based on insights from these experiments, in Section \ref{sec:adaptive}, we develop an adaptive version of the fault-tolerant algorithm which can dramatically reduce the number of false positive detections.

\section{Background}\label{sec:background}
\subsection{Silent errors}

Silent errors, also known as ``soft faults'' or ``silent data corruption'' (SDC), are faults that cause a change in some floating-point number without any apparent indication of a problem \cite{MeurantCG}. 
In this work, we assume that silent errors are \textit{transient}; that is, if an input to a computation is altered by a silent error, the result of the computation will be affected, but the error will not be persistent in the input after the computation. 
This assumption is common in practice, since inputs to a computation are in transient memory (cache); the same model is used, e.g., in \cite{ElliottIEEE} and \cite{CGSoft}.

We also assume that silent errors occur in the form of bit flips, since they are commonly studied \cite{ElliottIEEE} and easy to model.
Note that the impact of a bit flip may vary. A bit flip in the end of the mantissa of a floating-point number may have only negligible effects, whereas a flip of some dominant bit in the exponent can destroy the entire computation. 
It is possible that the prevalence of bit flips could increase in the future; as mentioned, as supercomputers become more complex and the number of their parts increases, the risk of a hardware failure increases as well \cite{MeurantCG}. Moreover, the relaxation of hardware correctness could be utilized as a way to save energy, and the decrease of transistor feature sizes makes individual components more prone to failure \cite{ElliottIEEE}.

\par \label{0110}
The most straightforward approaches for detecting silent errors are \textit{double modular redundancy} (DMR) and \textit{triple modular redundancy} (TMR)  \cite{CGSoft}. Here one performs the same computation multiple times, either consecutively on the same hardware unit or simultaneously on different hardware units, and then checks whether the results are the same \cite{CGSoft}. 
The crucial problem of redundancy approaches is their cost, as they require either multiple computational units or twice/thrice the time. This is especially limiting for large-scale parallel computers because of their energy consumption \cite{ElliottIEEE}.
\par
Another possible approach, which we use in this article, is to use information about the numerical method to derive a set of detection criteria \cite{MeurantCG}. This approach is called \textit{algorithm-based fault tolerance} (ABFT) \cite{MeurantCG}. 
Even though ABFT methods do not require the amount of computational resources needed for the redundancy approaches, they still require some. They may also cause delayed convergence if the algorithm is modified to  correct the errors during the computation \cite{MeurantCG}.

\subsection{The Pipe-PR-CG algorithm}

The HS-CG algorithm \cite{OGCG} is stated in Algorithm \ref{HSCG}. From now on, vector variables are written in bold and matrices in bold uppercase. Here, $\mathbf{M}$ denotes a symmetric positive definite preconditioner. Using $\mathbf{M}$, we can implicitly solve the system $\mathbf{L^{-T}AL^{-1}y} = \mathbf{L^{-T}b}$, where $\mathbf{y} = \mathbf{Lx}$ and $\mathbf{L}$ is the Cholesky factor of $\mathbf{M}$, utilizing just solutions of subsystems with $\mathbf{M}$ during the run, e.g., $\mathbf{\tilde{r}_k} = \mathbf{M^{-1}}\mathbf{r_k}$ \cite{ChenCarson}. The symbol $\sim$ denotes extra variables introduced by inclusion of the preconditioner. Here and in the remainder of the article we use $\langle \cdot, \cdot\rangle$ to denote the inner product of two vectors. The statement of the INITIALIZE() procedure can be found in Appendix \ref{init}.

\begin{algorithm}[H]
	\caption{Hestenes and Stiefel Conjugate Gradient: HS-CG (preconditioned)}
	\label{HSCG}
	\begin{algorithmic}[1]
		\Procedure {HS-CG}{$\mathbf{A}, \mathbf{M}, \mathbf{b}, \mathbf{x_0}$}
		\State INITIALIZE()
		\For {$ k = 1,2,\dots $}
		\State $\mathbf{x_k} = \mathbf{x_{k-1}} + \alpha_{k-1}\mathbf{p_{k-1}},\quad\mathbf{r_k} = \mathbf{r_{k-1}} - \alpha_{k-1}\mathbf{s_{k-1}},\quad  \mathbf{\tilde{r}_k} = \mathbf{M^{-1}}\mathbf{r_k}$
		\State $\nu_k = \langle \mathbf{\tilde{r}_k}$, $\mathbf{r_k} \rangle,\quad \beta_k = \nu_k/\nu_{k-1}$ \label{HSCG:nu}
		\State $\mathbf{p_k} = \mathbf{\tilde{r}_k} + \beta_{k}\mathbf{p_{k-1}}$, $\mathbf{s_k} = \mathbf{Ap_k}$
		\State $\mu_k = \langle \mathbf{p_k}, \mathbf{s_k} \rangle,\quad\alpha_k = \nu_k/\mu_k$
		\EndFor
		\EndProcedure
	\end{algorithmic}
\end{algorithm}
In Algorithm \ref{HSCG}, the computation has to be done largely sequentially, since each step depends on variables computed in previous steps.
This causes a communication bottleneck due to the global reductions needed to compute inner products and/or the distributed matrix-vector multiplication \cite{ChenCarson}.

The high cost of communication in HS-CG has motivated the search for mathematically equivalent CG variants better suited for implementation on parallel machines. One approach, called \textit{communication-hiding} variants (e.g., the variant of Chronopoulos and Gear \cite{ChGCG}), aim to reduce the number of synchronization points to one. 
One can go a step further, and ``pipeline'' the inner products and matrix operations to occur simultaneously, so that the global synchronization points no longer cause a bottleneck, as in the variant of Ghysels and Vanroose \cite{GVCG}. 
However, the rearrangement of computations and introduction of auxiliary recurrences necessary to achieve this can have negative affects on the convergence delay and attainable accuracy \cite{Carsonetal, ChenCarson}.

A pipelined variant which aims to mitigate these negative numerical effects is the Pipe-PR-CG method of \cite{ChenCarson}, shown in Algorithm \ref{Pipe-PRCG}. Pipe-PR-CG still requires only one global synchronization point per iteration. This is accomplished by deriving a mathematically equivalent expression for the variable $\nu_k$ (line \ref{HSCG:nu} in Algorithm \ref{HSCG}) utilizing quantities already computed in the previous iteration. This way we avoid the first computation of the inner product $\langle \mathbf{\tilde{r}_k}$, $\mathbf{r_k} \rangle$. However, this change could lead to a dramatic loss of accuracy as the value of $\nu_k$ could become negative \cite{ChenCarson}, and therefore we employ the new expression just as a ``predicted'' value, which is then used instead of the original expression to compute the next iterates until it is ``recomputed'' using the original inner product formulation. This idea was proposed by Meurant in \cite{MeurantCray} to stabilize the algorithm while also retaining the potential for parallelism. Although even this alteration might introduce some instability, it allows us to perform the iteration more efficiently on distributed memory machines since all the inner products occur at the same time \cite{ChenCarson}. Moreover, the maximal attainable accuracy is similar as for the original HS-CG. This perhaps surprising result is analyzed in \cite{ChenCarson}; we refer the reader to this work for a full derivation of the method. 

We focus on the Pipe-PR-CG method here for two reasons. First, it is the current state-of-the-art pipelined Krylov subspace method in terms of potential for parallel performance and numerical stability. Second, the predict-and-recompute aspect of the method provides ample opportunities for developing silent error detection approaches based on finite precision bounds. In short, the ``predicted'' and ``recomputed'' variables, which should be the same in exact arithmetic, provide quantities which can be easily compared in order to detect silent errors.

\begin{algorithm}[H]
	\caption{Pipelined Predict-and-Recompute Conjugate Gradient: \newline Pipe-PR-CG (preconditioned)}
	\label{Pipe-PRCG}
	\begin{algorithmic}[1]
		\Procedure {Pipe-PR-CG}{$\mathbf{A}, \mathbf{M}, \mathbf{b}, \mathbf{x_0}$}
		\State INITIALIZE()
		\For {$ k = 1,2,\dots $}
		\State $\mathbf{x_k} = \mathbf{x_{k-1}} + \alpha_{k-1}\mathbf{p_{k-1}}$
		\State $\mathbf{r_k} = \mathbf{r_{k-1}} - \alpha_{k-1}\mathbf{s_{k-1}},\quad  \mathbf{\tilde{r}_k} = \mathbf{\tilde{r}_{k-1}} - \alpha_{k-1}\mathbf{\tilde{s}_{k-1}}$
		\State $\mathbf{w'_k} = \mathbf{w_{k-1}} - \alpha_{k-1}\mathbf{u_{k-1}},\quad  \mathbf{\tilde{w}'_k} = \mathbf{\tilde{w}_{k-1}} - \alpha_{k-1}\mathbf{\tilde{u}_{k-1}}$
		\State $\nu'_k = \nu_{k-1} - \alpha_{k-1} \sigma_{k-1}-\alpha_{k-1} \phi_{k-1}+\alpha^2_{k-1}\gamma_{k-1},\quad\beta_k = \nu'_k/\nu_{k-1}$
		\State $\mathbf{p_k} = \mathbf{\tilde{r}_k} + \beta_{k}\mathbf{p_{k-1}},\quad\mathbf{s_k} = \mathbf{w'_k}+\beta_{k}\mathbf{s_{k-1}},\quad\mathbf{\tilde{s}_k} = \mathbf{\tilde{w}'_k}+\beta_{k}\mathbf{\tilde{s}_{k-1}}$
		\State $\mathbf{u_k} = \mathbf{A \tilde{s}_k},\quad \mathbf{\tilde{u}_k} = \mathbf{M^{-1}u_k},\quad\mathbf{w_k} = \mathbf{A \tilde{r}_k},\quad \mathbf{\tilde{w}_k} = \mathbf{M^{-1}w_k}$
		\State $\mu_k = \langle \mathbf{p_k}, \mathbf{s_k} \rangle,\, \sigma_k = \langle \mathbf{\tilde{r}_k}, \mathbf{s_k} \rangle,\, \phi_k = \langle \mathbf{\tilde{s}_k}, \mathbf{r_k} \rangle,\, \gamma_k = \langle \mathbf{\tilde{s}_k}, \mathbf{s_k} \rangle,\, \nu_k = \langle \mathbf{\tilde{r}_k}, \mathbf{r_k} \rangle$ \label{Pipe-PRCG-nu}
		\State $\alpha_k = \nu_k/\mu_k$
		\EndFor
		\EndProcedure
	\end{algorithmic}
\end{algorithm}

\section{The Effect of Silent Errors}\label{sec:effect}
In this section we investigate the sensitivity of Pipe-PR-CG (Algorithm \ref{Pipe-PRCG}) to silent errors. 
In the following text, the symbol $||\cdot||$ denotes the 2-norm. The condition number $\kappa(\mathbf{A})$ of a matrix $\mathbf{A}$ is defined as $\kappa(\mathbf{A}) = ||\mathbf{A}||\cdot||\mathbf{A}^{-1}||$.

Algorithm \ref{Pipe-PRCG} was implemented in Python (version 3.10.4) and experiments were performed on a computer with an 11th Generation Intel\textsuperscript{\textregistered} Core\textsuperscript{\texttrademark} i7-1185G7 processor and 16 GB of RAM running on 64-bit Windows 10 Pro operating system.
In all runs, the initial guess $\mathbf{x_0}$ was a vector of all zeros and the right-hand side $\mathbf{b}$ was such that
the vector of all ones $\mathbf{e} = (1,1,\cdots,1,1)^T$ was the exact solution of the system, i.e., $\mathbf{b=Ae}$. { This choice of the exact solution $\mathbf{e}$ may, depending on the matrix $\mathbf{A}$, yield some special right-hand sides. Nonetheless, setting $\mathbf{b}$ this way allows us to present a simple demonstration of sensitivity of Pipe-PR-CG to bit flips while having some degree of variety in the right-hand side.} No preconditioners were used, i.e., $\mathbf{M = I}$ and the computation of the variables with the tilde symbol (e.g., $\mathbf{\tilde{r}_k}$) is omitted since they are the same as their unpreconditioned counterparts.
The stopping criterion used was  $||\mathbf{r_k}||/||\mathbf{b}|| \le \epsilon_\text{tol}$ for the computation to conclude earlier than at the maximal allowed number of iterations, with $\epsilon_\text{tol}$ being $1\mathrm{e}{-10}$. 
Here we explicitly compute $\|\mathbf{r_k}\|$ to check for convergence rather than using the quantity $\nu_k$ (which in exact arithmetic gives the norm of $\mathbf{r_k}$), in order to have ``illustrative'' results for the case where a bit flip occurs in the variable $\nu_k$. 
\par
Additionally, the code was written so that each time an overflow warning occurs, an error is raised instead. This was done to suppress situations where the Python compiler does not terminate the computation right away, but assigns the result to infinity instead. There were also pure overflow errors, which stopped the computation immediately. We count all these cases as ``did not converge''. In \cite{CGSoft}, the same categorization of overflows as ``non-convergent'' was used. This is the reason why there were some non-convergent cases for the variable $\mathbf{x_k}$, even though it does not influence
any other variables, and therefore we would expect the runs with silent errors in it
to be always ``convergent''. 

The injection of silent errors into variables was implemented using the Python module \textit{bitstring} (version 4.1) \cite{bitstring}. Time-wise, the flips always occur after the new value of a variable is computed; e.g., we first compute $\alpha_k = \nu_k / \mu_k$, and then insert a bit flip into $\alpha_k$.
For each matrix, a run with no bit flips was performed to determine the number of
iterations $\varphi$ needed to converge. A run ``tainted'' by a silent error was then deemed
as ``converged'' if it reached the stopping criteria within $1.5\varphi$ iterations. The same approach for determining convergence was used in the article \cite{CGSoft}.

The experiment was performed for each of the 14 variables in Pipe-PR-CG ($\mathbf{x_k}$, $\mathbf{r_k}$, $\mathbf{w'_k}$,  $\nu'_k$, $\beta_k$, $\mathbf{p_k}$, $\mathbf{s_k}$, $\mathbf{u_k}$, $\mathbf{w_k}$, $\mu_k$,  $\sigma_k$, $\gamma_k$, $\nu_k$, $\alpha_k$). {We note that in the unpreconditioned case, we can take $\phi_k=\mu_k$. }
We tested 3 different variants of when the bit flip occurred: $0.3\varphi$, $0.6\varphi$, and
$0.9\varphi$ iterations. This was performed for all 64 bits. In case of scalar variables, one run was
performed for each matrix from the dataset, bit number, and flip iteration. For vector variables,  the index in which the bit flip occurs was chosen randomly. There were 20 trials
for each bit, flip iteration, and matrix, so that the randomness in the index choice
could be included.

We tested a number of SPD matrices from the SuiteSparse Matrix Collection \cite{SuiteSparse,SuiteSparseArticle}, listed in Table \ref{tab:matrices}. The matrices used were selected to represent various sizes, condition numbers, singular value distributions, structures, as well as problem sources. 

\begin{table}[]
\centering
\begin{tabular}{|c|c|c|}
\hline
name & $n$ & $\kappa(A)$ \\ \hline
\textit{1138\_bus} & 1,138 & $8.6e+06$ \\ \hline
\textit{bcsstm07} & 420 & $7.6e+03$ \\ \hline
\textit{bundle1} & 10,581 & $1.0e+03$ \\ \hline
\textit{wathen120} & 36,441 & $2.6e+03$ \\ \hline
\textit{bcsstk05} & 153 & $1.4e+04$ \\ \hline
\textit{gr\_30\_30} & 900 & $1.9e+02$ \\ \hline
\textit{nos7} & 729 & $2.4e+09$ \\ \hline
\textit{crystm01} & 4,875 & $2.3e+02$ \\ \hline
\textit{aft01} & 8,205 & $4.4e+18$ \\ \hline
\end{tabular}
\caption{Matrices used in experiments}
\label{tab:matrices}
\end{table}

\begin{figure}[!htb]
    \begin{minipage}[t]{0.48\textwidth}
        \centering
        \includegraphics[clip, trim=0cm 0cm 0cm 0cm, width=1.1\textwidth]{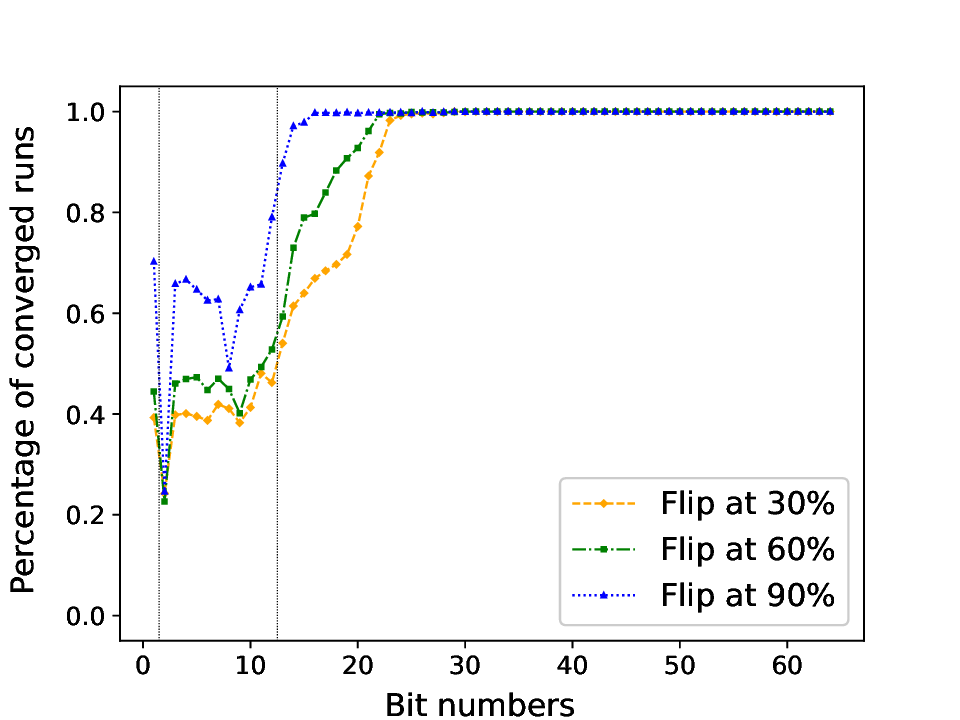}
        \caption{Bit flip sensitivity: averages from all variables.}
        \label{Fig1}
    \end{minipage}
    \hspace{5pt}
    \begin{minipage}[t]{0.48\textwidth}
        \centering
        \includegraphics[clip, trim=0cm 0cm 0cm 0cm, width=1.1\textwidth]{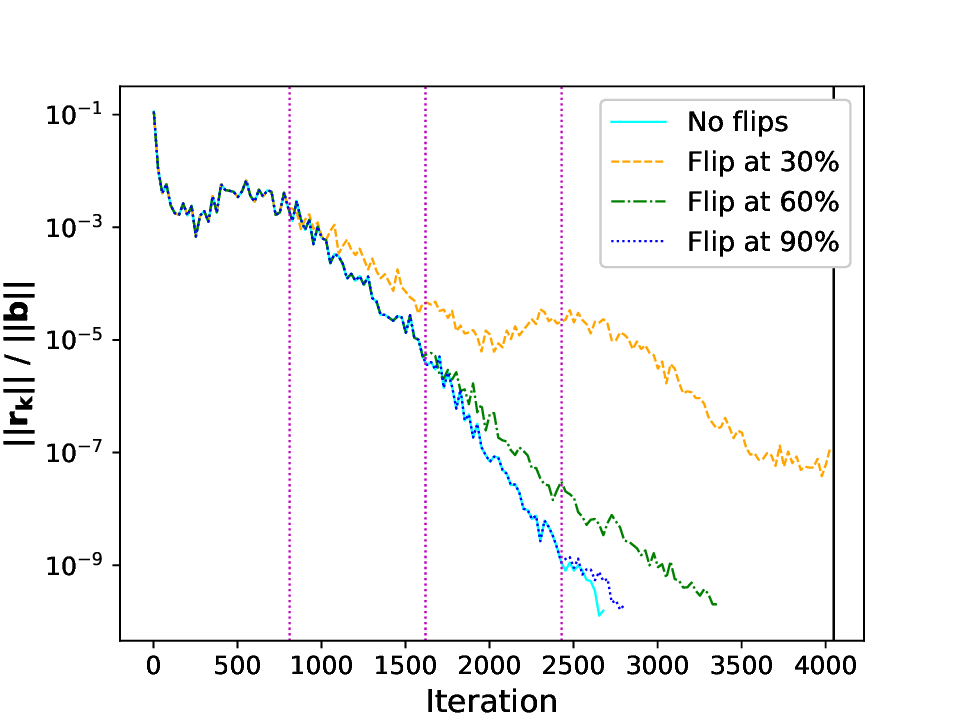}
        \caption{Convergence of the relative residual $\mathbf{||r_k||/||b||}$ when the 15th bit of $\beta_k$ is flipped for matrix \textit{1138\_bus}. Purple dotted lines denote flip iterations.}
        \label{Fig2}
    \end{minipage}
\end{figure}

The output of the experiment is a graph depicting what percentage of runs are ``convergent'' for each of the 64 bits. 
The results are presented as averages over all variables in Figure \ref{Fig1}. The averages are calculated with each variable having the same weight. The graph includes thin vertical lines, which separate the sign (1 bit), exponent (11 bits), and
mantissa (52 bits). Note that here and in the remainder of the article the bits are numbered from 1 to 64 using big-endian ordering (i.e., the leftmost bit is 1).

In Figure \ref{Fig1}, spikes in the curves might be caused by the fact that bit flips from 0 to 1 and from 1 to 0 are not equally significant \cite{CGSoft}.
A non-convergent spike can be observed for the second bit, but this might be expected as it is the most significant bit in terms of the absolute value of a number. Therefore, it is most likely to cause an overflow error that inflates the non-convergent cases. The small drops for the 8th bit for flip at $0.9\varphi$ and the 9th bit for flips at $0.3\varphi$ and $0.6\varphi$ might be caused by these bits being significant for some value range our variables often fall into.

{ It is well-known that the earlier a bit flip occurs, the greater its effect on the overall convergence; see, e.g., \cite{MeurantCG}.} We examine this effect more closely in Figure \ref{Fig2}, which contains convergence curves of the relative residual $||\mathbf{r_k}||/||\mathbf{b}||$ in Pipe-PR-CG for all three time-wise flip options when the 15th bit of $\beta_k$ is flipped for matrix \textit{1138\_bus}. Note that to facilitate legibility of the plot, only every 25th iteration is plotted for the convergence curves. Dotted purple lines and the solid black line denote when the flips occurred and where the $1.5\varphi$ termination point is, respectively. We can indeed see that in this case the computation was more heavily influenced by an earlier flip. On the other hand, flipping the 15th bit at a point when the method had almost converged did not have a significant effect on the number of extra iterations necessary. 
\par
In general, it seems that silent errors in bits numbered around 25 and higher have little influence on convergence delay. This might be expected due to the somehow decreasing ``significance'' of bits as we proceed to those with higher index. A special qualitatively different case is the first bit - the sign - as it, unlike the other bits, does not influence the absolute value of the number.

\begin{figure}[!htb]
	\centering
	\begin{subfigure}{.5\textwidth}
		\centering
		\includegraphics[clip, trim=0cm 0cm 1cm 0.5cm,width=\linewidth]{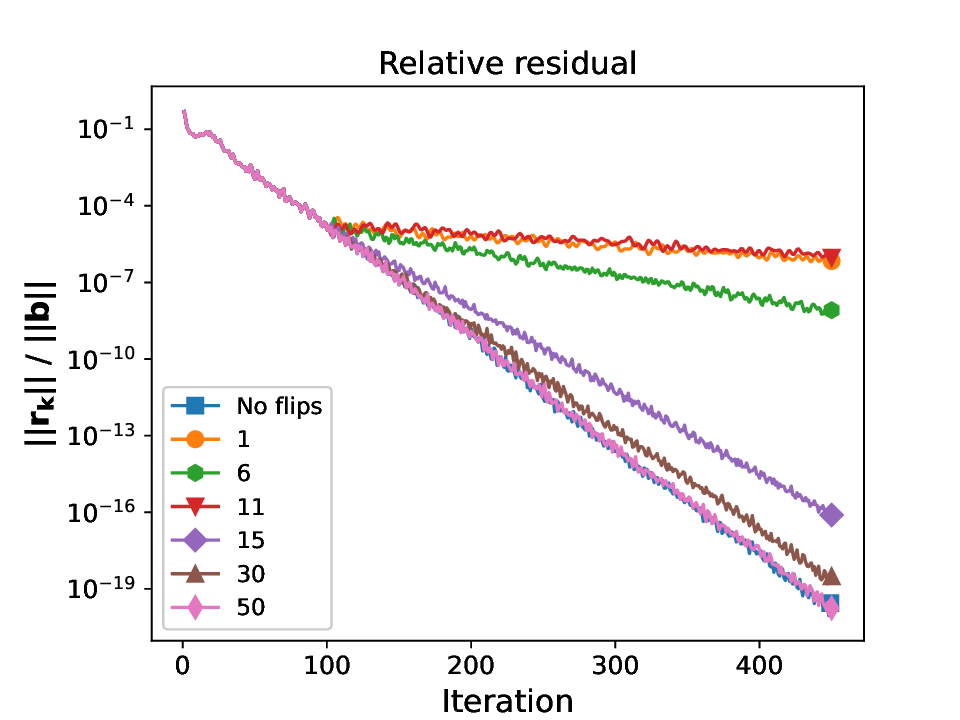}
	\end{subfigure}%
	\begin{subfigure}{.5\textwidth}
		\centering
		\includegraphics[clip, trim=0cm 0cm 1cm 0.5cm,width=\linewidth]{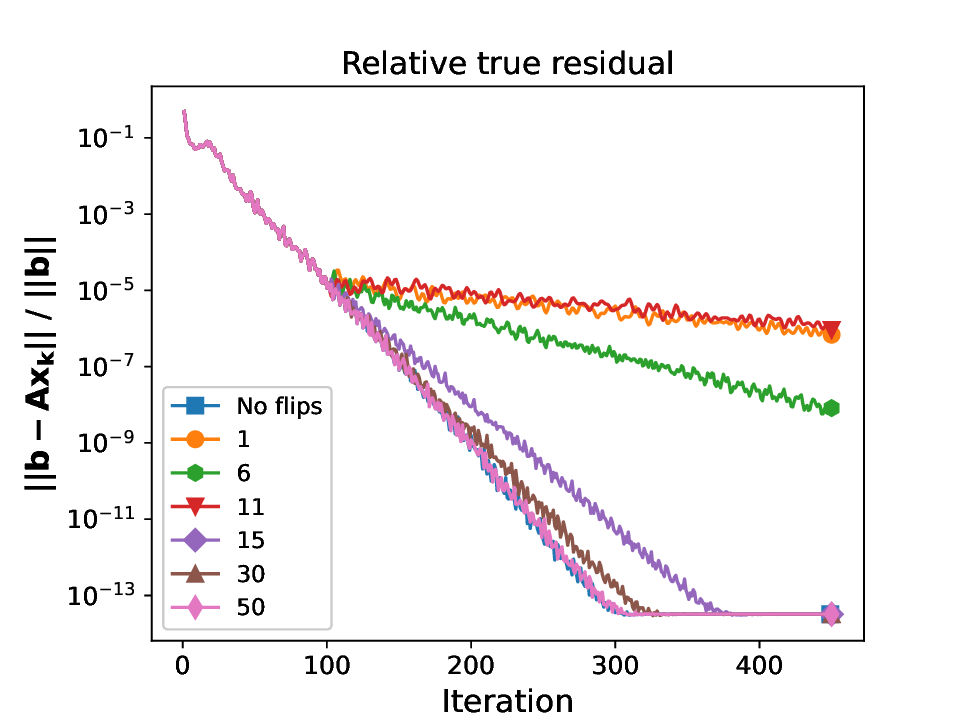}
	\end{subfigure}
	\caption{Residual convergence curves when various bits of $\sigma_k$ are flipped in the 100th iteration for the matrix \textit{bundle1}.}
	\label{Fig3}
\end{figure}

The decreasing effect of a bit flip as the bit number rises is illustrated in Figure \ref{Fig3}, which, as an exception in this section, has a fixed number of iterations for each run in order to illustrate the behavior more clearly. We can see that for our data, in terms of both the relative residual (left plot) and the relative true residual (right plot), the computation is impacted more significantly by flips in the sign and exponent bits. The fact that the flip in the 11th bit is more influential in this case than the one in the 6th bit may again be caused by whether it is a ``from 0 to 1'' or ``from 1 to 0'' flip.
\par
Thus as expected, bit flips influence Pipe-PR-CG more significantly when they occur early in the computation or when they are in bits which are either the sign or have a serious impact on the absolute value of the altered number.  

\subsection{Other effects and error detection}
As we have seen in the previous section, silent errors can have significant influence on the convergence of Pipe-PR-CG. Naturally, it can be surmised that other aspects of the procedure might be affected as well. Consequently, some of the effects could be utilized for our ultimate goal - the detection of silent errors. This is the idea of the \textit{algorithm-based fault tolerance methods}. The fundamental concept of this approach is to derive some criteria of silent error detection from the theoretical or practical knowledge we possess of the algorithm \cite{MeurantCG}. In our case, we will try to utilize the predict-and-recompute principle which allows us to hide some of the communication. 
\par
In Pipe-PR-CG, there are two variables whose value is first predicted using an alternative relation and then recomputed,  $\nu_k$ and $\mathbf{w_k}$. We will investigate the ``gaps'' between their predicted and recomputed versions, i.e., $|\nu_k - \nu_k'|$ and $||\mathbf{w_k} - \mathbf{w_k'}||$, which should be zero in exact arithmetic. We begin with the \textit{$\nu$-gap}. 
\par
\begin{figure}[!htb]
	\makebox[\textwidth][c]{\includegraphics[clip, trim=0cm 0cm 0cm 1cm, width=.6\textwidth]{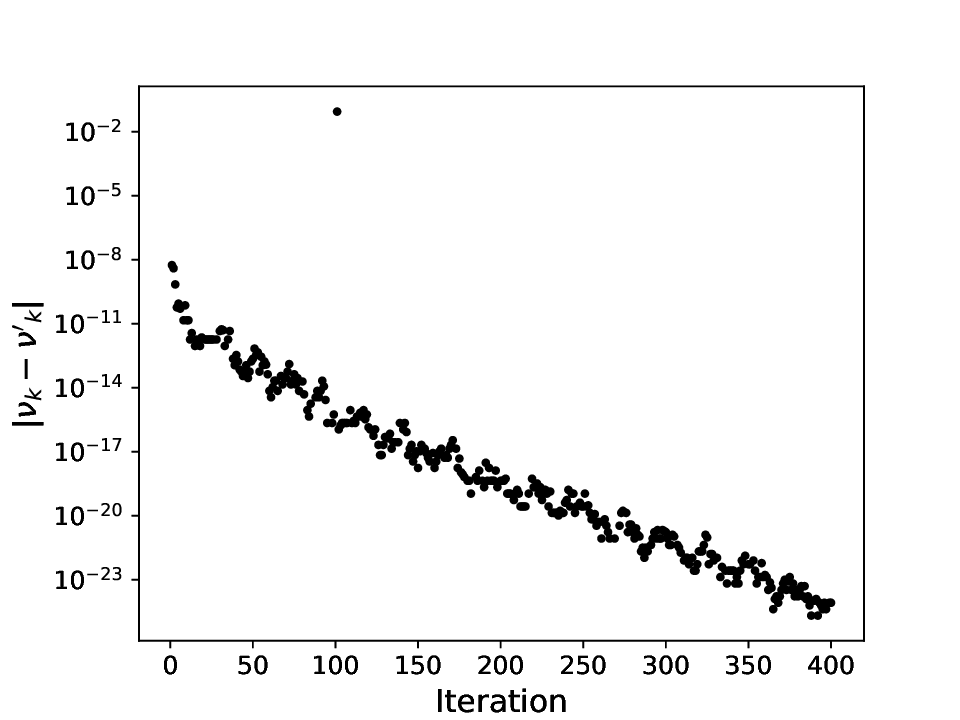}}
	\caption{The $\nu$-gap ($|\nu_k - \nu_k'|$) when the 15th bit of $\gamma_k$ is flipped in the 100th iteration for the matrix \textit{bcsstm07}}
	\label{Fig4}
\end{figure}

Figure \ref{Fig4} depicts the size of the $\nu$-gap when the 15th bit of $\gamma_k$ is flipped in the 100th iteration for the matrix \textit{bcsstm07}. The computation is without preconditioning, with initial guess $\mathbf{x_0}$ being a vector of all zeros, and the right-hand side $\mathbf{b}$ is chosen as before. 
As can be observed, there is a significant outlier among the values in the 101st iteration, i.e., in the very next iteration after the bit flip. This is a promising result which indicates that the $\nu$-gap might be useful for silent error detection and motivates further exploration. The next section studies the possibility of using several variable ``gaps'' to this end. 

\section{Detecting Silent Errors}\label{sec:detection}
In this section, we derive relations which can be effectively used to detect silent errors in Pipe-PR-CG (Algorithm \ref{Pipe-PRCG}). We present several so-called ``gaps'', and then employ rounding error analysis to obtain expressions which bound these gaps from above. Subsequently, numerical experiments are performed for each of these gap-bound pairs to judge how effectively they detect silent errors. The idea is that a violation of the bound is a potential indicator of a bit flip having occurred. 
\par

{In the following, we will present rounding error bounds for both the unpreconditioned and the preconditioned Pipe-PR-CG algorithms. For our experiments, however, we will test the unpreconditioned case for simplicity. }

We will use a standard model of arithmetic with floating-point numbers. Let the symbol $\circ$ denote one of the operations $\{+,-,\times,\div\}$, $\epsilon$ the machine precision, $a$ and $b$ arbitrary feasible real numbers, and $\text{fp}(\cdot)$ which operation is performed in finite precision. Then, it holds that
\begin{equation} \label{41}
	|\text{fp}(a \circ b) - a \circ b| \le \epsilon \,|a \circ b|.
\end{equation}
Within this framework, it is possible to derive bounds on some of the standard vector operations. Letting $\mathbf{x},\mathbf{y} \in \mathbb{R}^n$ and $a \in \mathbb{R}$, we have that
\begin{align} 
	||\text{fp}(\mathbf{x} + a\mathbf{y}) - (\mathbf{x} + a\mathbf{y})|| &\le \epsilon \,(||\mathbf{x}|| + 2|a|\,||\mathbf{y}||),\label{42_1}\\
	||\text{fp}(\langle \mathbf{x}, \mathbf{y} \rangle) - \langle \mathbf{x}, \mathbf{y} \rangle|| &\le \epsilon \, n \, ||\mathbf{x}|| \,||\mathbf{y}||,\label{42_2}\\
	||\text{fp}(\mathbf{Ax}) - \mathbf{Ax}|| &\le \epsilon \, c \, ||\mathbf{A}|| \,||\mathbf{x}||,\label{42_3}
\end{align}
where $c$ is a constant depending on specific properties of the matrix $\mathbf{A}$. For instance, it is frequently taken as $c=mn^{1/2}$, where $m$ is the maximum number of nonzeros over the rows of $\mathbf{A}$. 
\par
Here, we use $\delta$ to denote a round-off error introduced in a calculation, i.e., the difference between the actual computed result and the exact expression for the variable denoted in the subscript of $\delta$. Taking for instance $\delta_{\sigma_k}$ as an example, it holds that $\delta_{\sigma_k} = \sigma_k - \langle \mathbf{r_k}, \mathbf{s_k} \rangle \le \epsilon \, n \, ||\mathbf{r_k}|| \,||\mathbf{s_k}||$.

\subsection{$\nu$-gap}
We first investigate the $\nu$-gap mentioned in the previous section, defined as the size of the difference between the predicted value $\nu'_k$ and the recomputed value $\nu_k$ in Pipe-PR-CG, i.e., $|\nu_k - \nu'_k|$. Let us recall how these variables are defined. It holds that
\[\nu_k = \langle \mathbf{r_k}, \mathbf{r_k} \rangle, \quad\text{and}\quad \nu'_k = \nu_{k-1} - \alpha_{k-1} \sigma_{k-1}- \alpha_{k-1} \phi_{k-1}+\alpha^2_{k-1}\gamma_{k-1},\]
which are mathematically equivalent, i.e., they are equal in exact arithmetic.
As we observed in Figure \ref{Fig4}, the $\nu$-gap shows promising potential for silent error detection. However, an issue is how to determine when the value of the $\nu$-gap signals a potential silent error occurrence. One problem is that, as can be seen in Figure \ref{Fig4}, the $\nu$-gap can fluctuate. Moreover, even if a silent error influences the $\nu$-gap, nothing guarantees there will always be such a distinct outlier value as in the aforementioned graph. Therefore, we aim to derive a bound on the $\nu$-gap, and use this bound to determine if a silent error occurred by checking whether the computed $\nu$-gap exceeds the bound.

A bound for the $\nu$-gap {for the unpreconditioned case} was previously derived in \cite{ChenCarson}. By using slightly different algebraic manipulations than those in \cite{ChenCarson}, we can derive the tighter bound 
\[|\nu_k - \nu_k'|\equiv\Delta_{\nu'_k} \lesssim \epsilon \, (21+6n)(||\mathbf{r_{k-1}}||^2 + ||\mathbf{r_{k}}||^2),\]
which has slightly smaller constants than the bound given in \cite{ChenCarson}. Note that the $\lesssim$ indicates that we have dropped terms of order $O(\epsilon^2)$ and higher. 

{
By applying a similar rounding error analysis to the preconditioned variant of Pipe-PR-CG, we can derive a bound for the $\nu$-gap as:
\[
    \Delta_{\nu_k'} \lesssim \epsilon \left( \frac{21+6n}{2} \right) \left(\|\rrkm\|^2 + \|\trkm\|^2 + \|\rrk\|^2 + \|\trk\|^2 \right).
\]
}
{ Note that the norms $\|\rrk\|$ and $\|\trk\|$ can be computed at the same time as the other inner products in line \ref{Pipe-PRCG-nu} of Algorithm \ref{Pipe-PRCG} without compromising the potential for pipelining.} 

\begin{figure}[h]
	\makebox[\textwidth][c]{\includegraphics[clip, trim=0cm 0cm 0cm 1cm, width=.7\textwidth]{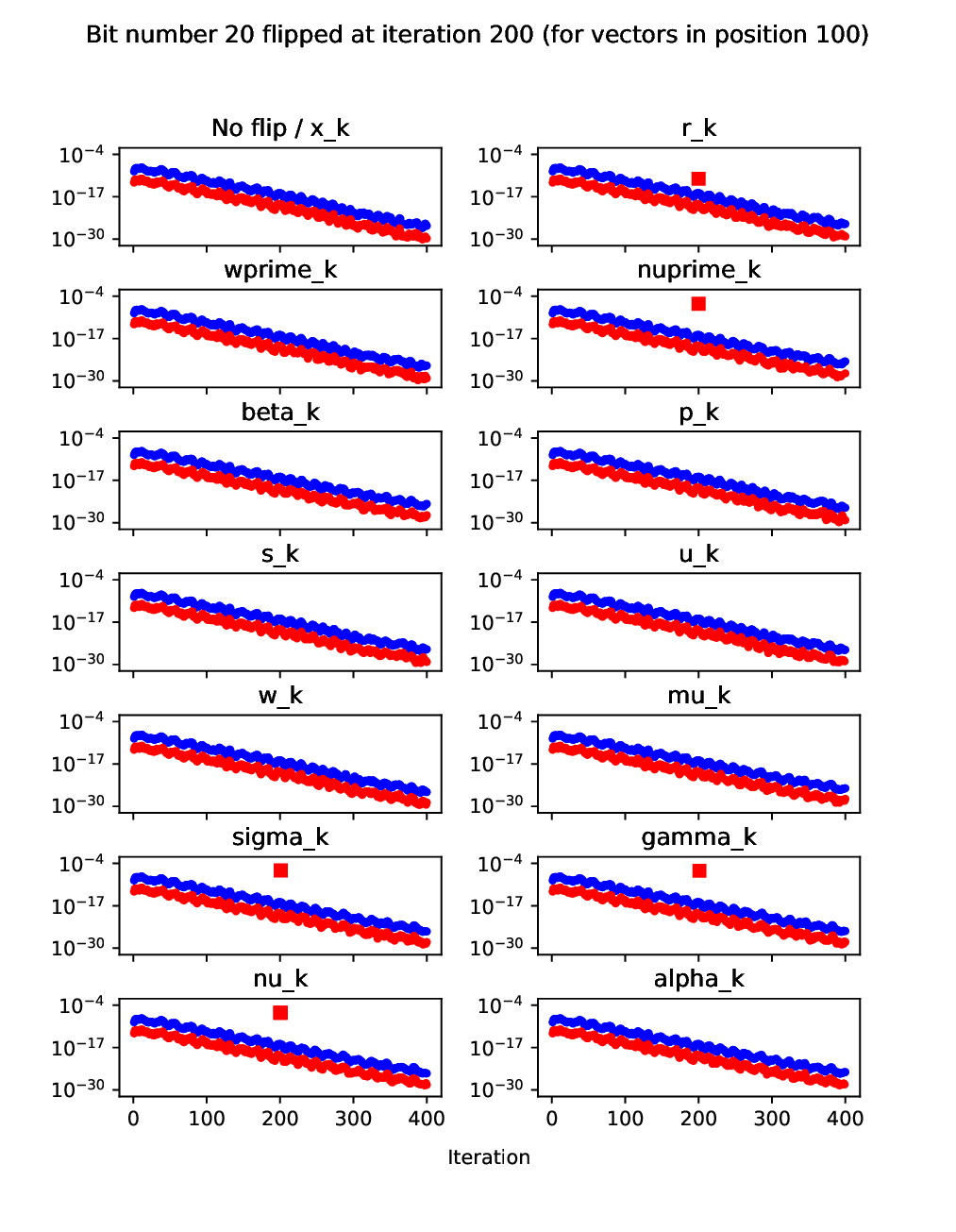}}
	\caption{$\nu$-gap (red) and $\nu$-gap bound (blue) graph, matrix \textit{bcsstm07}}
	\label{Fig5}
\end{figure}

We now investigate whether this bound is violated if a bit flip occurs. This is illustrated in Figure \ref{Fig5}, which depicts the behavior of the $\nu$-gap, $|\nu_k - \nu_k'|$, and the derived bound on the $\nu$-gap, $\epsilon \, (21+6n)(||\mathbf{r_{k-1}}||^2 + ||\mathbf{r_{k}}||^2)$, when the 20th bit is flipped in each variable in the 200th iteration of Pipe-PR-CG for the matrix \textit{bcsstm07} (see Table \ref{tab:matrices}). The right-hand side $\mathbf{b}$ was a vector of all ones, i.e., $\mathbf{b} = \mathbf{e}$. There is one common subplot for the variable $\mathbf{x_k}$ and the case of no bit flip occurring during the run, since the resulting graph is the same since flips in $\mathbf{x_k}$ have no effect on the values of other variables. 
The norms $||\mathbf{r_{k-1}}||^2$ and $||\mathbf{r_k}||^2$ were computed using the already calculated $\nu_k = \langle \mathbf{r_k}, \mathbf{r_k} \rangle$ and $\nu_{k-1}=\langle \mathbf{r_{k-1}}, \mathbf{r_{k-1}} \rangle$.
For vector variables ($\mathbf{x_k}$, $\mathbf{r_k}$, $\mathbf{w'_k}$, $\mathbf{p_k}$, $\mathbf{s_k}$, $\mathbf{u_k}$, $\mathbf{w_k}$) the bit was always flipped in the 100th position of the vector.

When the $\nu$-gap exceeds the $\nu$-gap bound, its marker symbol is changed to a square for visibility.
The $\nu$-gap can at times be zero. However, this is not displayed in the figures for simplicity.  Finally, note that many more experimental runs with different system data {(including different right-hand sides)} were performed in order to judge the behavior of the $\nu$-gap and the $\nu$-gap bound; the plot presented here is a representative sample. This holds for all the silent error detection criteria introduced in this section; additional experiments for a different matrix can be found in Appendix \ref{AppendixFigures}.

Generally, it can be concluded that the $\nu$-gap detection works well when the silent error occurs in $\nu'_k,\sigma_k,\gamma_k$, or $\nu_k$. It can also seemingly detect errors in the residual $\mathbf{r_k}$. However, this does not hold if the flip occurs in the first (sign) bit. This is only logical, as the sign bit is irrelevant for the value of $\nu_k = \langle \mathbf{r_k}, \mathbf{r_k} \rangle$, and therefore the $\nu$-gap is not influenced by flips in it.
\par
The bound is violated even for flips in bits of higher number. During the sensitivity experiments in the previous section, it was discovered that flips in bits of number 25 and higher usually do not destroy convergence, as illustrated by Figure \ref{Fig1}. This means that once we employ this criterion in practice, it may raise an alarm even in the case of bit flips which have a negligible effect on convergence. Also noteworthy is that for $\mathbf{r_k}$ and $\nu'_k$, the bound is violated at the flip iteration, whereas for $\sigma_k$ and $\gamma_k$, the violation happens one iteration later. In the case of $\nu_k$, the bound is violated both in the flip iteration and the following iteration. Monitoring when our criteria raise an alarm that a silent error has likely occurred is important for its correction. The idea for this correction is that we keep variables from a number of previous iterations or make some checkpoints, and if the flip is detected immediately we can roll back to a state which should not yet be influenced by the error. 
\par
In conclusion, the $\nu$-gap criterion seems to be able to reliably detect flips in $\nu'_k,\sigma_k,\gamma_k$, and $\nu_k$, as well as in non-sign bits for $\mathbf{r_k}$. For other variables a different detection method must be used.

\subsection{$w$-gap}
We can also investigate the $\mathbf{w}$-gap, i.e., the size of the difference between the predicted value $\mathbf{w'_k}$ and the recomputed value $\mathbf{w_k}$, which again are equal in exact arithmetic. Note that there is no existing bound on this quantity in the literature. 
\par
We now seek to derive a bound for the $\mathbf{w}$-gap. With the symbols $\delta$ once again denoting the rounding errors, it holds that 
\begin{align} \label{413}
	& \mathbf{r_k} = \mathbf{r_{k-1}} - \alpha_{k-1}\mathbf{s_{k-1}} + \delta_{\mathbf{r_k}}, \quad \mathbf{w'_k} = \mathbf{w_{k-1}} - \alpha_{k-1}\mathbf{u_{k-1}} + \delta_{\mathbf{w'_k}}, \\
	& \mathbf{u_k} = \mathbf{A s_k} + \delta_{\mathbf{u_k}}, \quad
	\mathbf{w_k} = \mathbf{A r_k}+ \delta_{\mathbf{w_k}} \nonumber.
\end{align}
We use these relations to rewrite the expression for $\mathbf{w_k} - \mathbf{w'_k}$ as
\begin{align*}
	\mathbf{{w}_{k}} -  \mathbf{{w}'_{k}} & = \mathbf{Ar_k} + \delta_{\mathbf{w_k}} - (\mathbf{w_{k-1}} - \alpha_{k-1}\mathbf{{u}_{k-1}} + \delta_{\mathbf{w}'_k})\\
	& = \mathbf{Ar_{k-1}} - \alpha_{k-1}\mathbf{As_{k-1}} + \mathbf{A}{\delta_{\mathbf{r_k}}} + \delta_{\mathbf{w_k}} - (\mathbf{w_{k-1}} - \alpha_{k-1}\mathbf{{u}_{k-1}} + \delta_{\mathbf{w}'_k})\\
	& = (\mathbf{Ar_{k-1}} - \mathbf{w_{k-1}}) - \alpha_{k-1}(\mathbf{As_{k-1}} -\mathbf{{u}_{k-1}}) + \mathbf{A}{\delta_{\mathbf{r_k}}} + \delta_{\mathbf{w_k}} -  \delta_{\mathbf{w}'_k}\\
	& = - \delta_{\mathbf{{w}_{k-1}}} + \delta_{\mathbf{w_k}} + \alpha_{k-1}\delta_{\mathbf{{u}_{k-1}}} + \mathbf{A}{\delta_{\mathbf{r_k}}} -  \delta_{\mathbf{w}'_k}.
\end{align*}
Subsequently, we can take the norm of both sides to obtain the inequality
\begin{equation} \label{414}
	\Delta_{\mathbf{w}'_k} \equiv ||\mathbf{{w}_{k}} -  \mathbf{{w}'_{k}} || \le ||\delta_{\mathbf{{w}_{k-1}}}|| + ||\delta_{\mathbf{w_k}}|| + |\alpha_{k-1}|\,||\delta_{\mathbf{{u}_{k-1}}}|| + ||\mathbf{A}||\,||\delta_{\mathbf{r_k}}|| +  ||\delta_{\mathbf{w}'_k}||.
\end{equation}
We now aim to bound individual terms from \eqref{414}, starting with the first three of them, i.e., $\delta_{\mathbf{{w}_{k-1}}}$, $\delta_{\mathbf{{w}_{k}}}$, and $\delta_{\mathbf{{u}_{k-1}}}$. Using relations from \eqref{413} together with \eqref{42_3} yields
\begin{align}
	||\delta_{\mathbf{{w}_{k}}}|| &\le \epsilon \, c \, ||\mathbf{A}||\,||\mathbf{r_k}||,\label{415_1}\\
	||\delta_{\mathbf{{w}_{k-1}}}|| &\le \epsilon \, c \, ||\mathbf{A}||\,||\mathbf{r_{k-1}}||,\label{415_2}\\
	||\delta_{\mathbf{{u}_{k-1}}}|| &\le \epsilon \, c \, ||\mathbf{A}||\,||\mathbf{s_{k-1}}||\label{415_3},
\end{align}
where $c = mn^{1/2}$, with $m$ the maximum number of nonzero elements in any row of $\mathbf{A}$.
\par 

If we rewrite the expression for $\mathbf{r_k}$ from \eqref{413} we obtain that
\[\mathbf{s_{k-1}} = \frac{1}{\alpha_{k-1}}(\mathbf{r_{k-1}}-\mathbf{r_k} + \delta_{\mathbf{r_k}}),\]
from which we can bound $\mathbf{s_k}$ from above as
\begin{equation} \label{46}
	\mathbf{||s_{k-1}}|| \le \frac{1}{|\alpha_{k-1}|}(||\mathbf{r_{k-1}}||+||\mathbf{r_k}|| + ||\delta_{\mathbf{r_k}}||).
\end{equation}

{Now, let us start by bounding $\delta_{\mathbf{r_k}}$. Using \eqref{42_2}, \eqref{413}, and \eqref{46}, and dropping $O(\epsilon^2)$ terms, we have
\begin{align}
	||\delta_{\mathbf{r_{k}}}|| 
	&\lesssim 3\epsilon \,(||\mathbf{r_{k-1}}||+||\mathbf{r_k}||).  
\end{align} \label{47}

We can bound $\delta_{\mathbf{w'_{k}}}$ using  \eqref{42_2}, \eqref{413}, \eqref{415_2}, \eqref{415_3}, and \eqref{46}, and again dropping $O(\epsilon^2)$ term, to obtain 
\begin{align}
	||\delta_{\mathbf{w'_{k}}}|| 
	& \lesssim 3\epsilon \, ||\mathbf{A}|| \,(||\mathbf{r_{k-1}}||+||\mathbf{r_{k}}||). \label{419}
\end{align} 
With each term bounded, we can now substitute \eqref{415_1} -  \eqref{419} into  \eqref{414}, and then drop terms of order $O(\epsilon^2)$ to obtain 
\begin{align*}
	\Delta_{\mathbf{w}'_k} &\lesssim 2(c+3)\epsilon||\mathbf{A}||\, (||\mathbf{r_{k-1}}||+||\mathbf{r_{k}}||).
\end{align*}}
{ Using a similar derivation, for the preconditioned version of the algorithm, the bound becomes
\[
	\Delta_{\mathbf{w}'_k} \lesssim 2(c+3)\epsilon||\mathbf{A}||\, (||\mathbf{\tilde{r}_{k-1}}||+||\mathbf{\tilde{r}_{k}}||).
\]
As before, the norms $||\trk||$ can be computed in each iteration along with the other inner products in line \ref{Pipe-PRCG-nu} of Algorithm \ref{Pipe-PRCG} without destroying the pipelining. 
}

Unlike the bound on the $\nu$-gap, these bounds contain two terms, $c$ and $||\mathbf{A}||$, which depend on properties of $\mathbf{A}$, and which ideally should be known prior to the computation if we wish to utilize this bound for silent error detection. 
We note that because these rounding error bounds are obtained via a worst-case analysis, they will often be very loose bounds in practice, and thus a rough estimate of $||\mathbf{A}||$ will likely suffice. If the user has no prior knowledge about $||\mathbf{A}||$, an inexpensive estimate could be obtained via sampling or randomized methods; see, e.g., \cite{martinsson2020randomized}.
Alternatively, for the preconditioned algorithm, we note that one could use the relation $\|\mathbf{A}\| = \|\mathbf{M}\mathbf{M^{-1}}\mathbf{A}\|\leq \|\mathbf{M}\| \|\mathbf{M^{-1}}\mathbf{A}\|$. A reasonable estimation of $||\mathbf{M^{-1}}\mathbf{A}||$ can be obtained inexpensively from a few iterations of Pipe-PR-CG itself; see the approach in \cite{Anorm}.
In the case that $\mathbf{M}$ consists of a decomposition into triangular factors, its norm can also easily be estimated.

\par
With the bound on the $\mathbf{w}$-gap derived, we can now test its efficacy for silent error detection. We again present a multi-graph as an example of the behavior, this time depicting the behavior of the $\mathbf{w}$-gap and the bound on the $\mathbf{w}$-gap when a bit flip occurs in each of the Pipe-PR-CG variables.
The norms $||\mathbf{r_k}||$ and $||\mathbf{r_{k-1}}||$ were computed by taking the square roots of  $|\nu_k|$.
The setup is the same as for the $\nu$-gap. 
\par
\begin{figure}[h]
	\makebox[\textwidth][c]{\includegraphics[clip, trim=0cm 0cm 0cm 1cm, width=.7\textwidth]{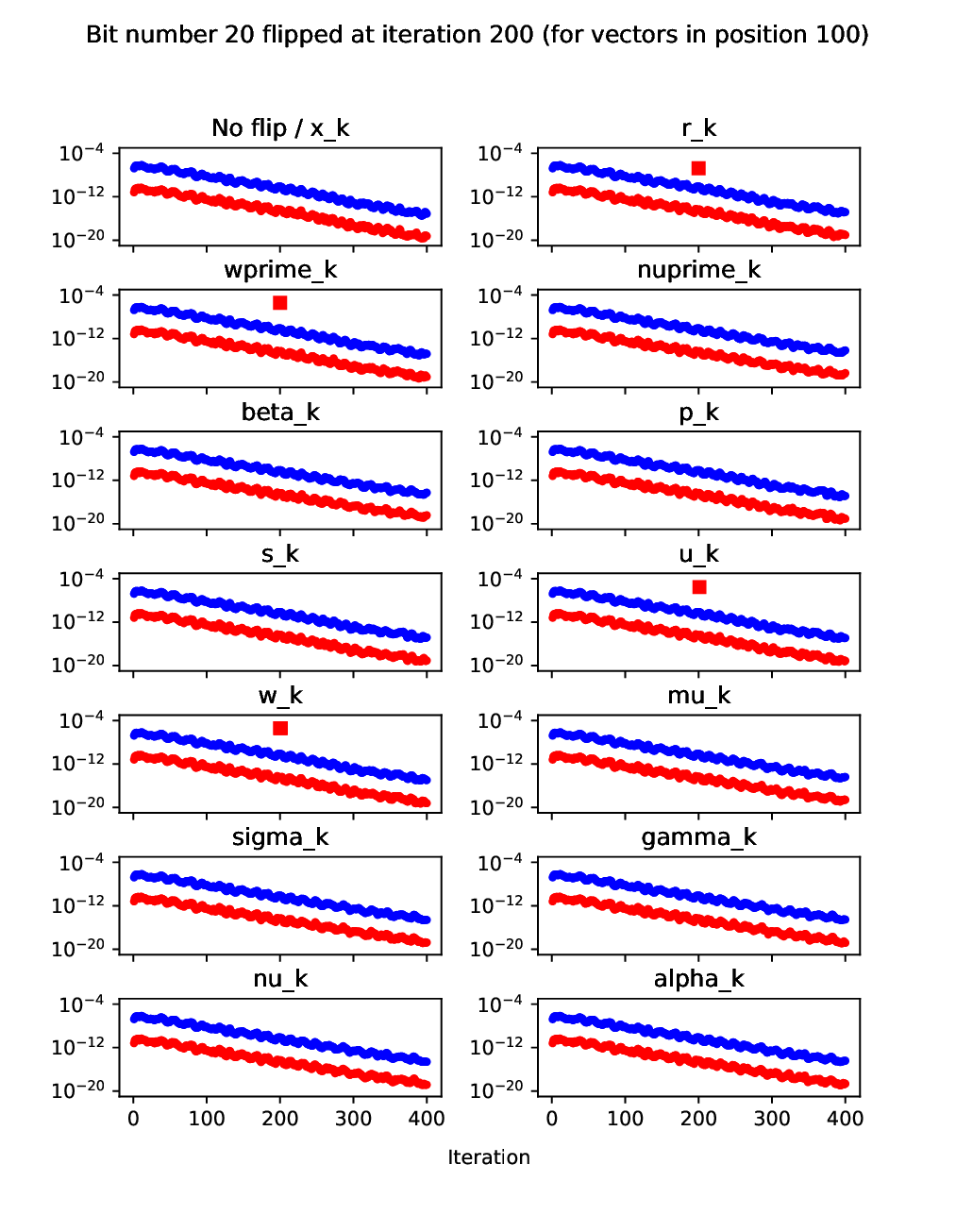}}
	\caption{$\mathbf{w}$-gap (red) and $\mathbf{w}$-gap bound (blue) graph, matrix \textit{bcsstm07}}
	\label{Fig6}
\end{figure}

In Figure \ref{Fig6}, we observe that the $\mathbf{w}$-gap bound seems to work for silent error detection in $\mathbf{r_k}$, $\mathbf{w'_k}$, $\mathbf{u_k}$, and $\mathbf{w_k}$. The $\mathbf{w}$-gap is also able to detect flips of the sign bit in the residual vector (data not shown). 
Once again, the bound is violated either at the flip iteration (for $\mathbf{r_k}$ and $\mathbf{w'_k}$), in the following iteration (for $\mathbf{u_k}$), or in both (for $\mathbf{w_k}$).

\subsection{$\mu$-gap}
Having investigated both the $\nu$-gap and the $\mathbf{w}$-gap, it is apparent that there are still some Pipe-PR-CG variables for which these two detection methods do not work. Therefore, it is necessary to derive another criterion which is not based on monitoring a difference between the predicted value and the recomputed value of some variable. Fortunately, there are other quantities in the Pipe-PR-CG algorithm which should be equal in exact arithmetic. These are $\mu_k$ and $\sigma_k$, defined as
\[\mu_k = \langle  \mathbf{p_k}, \mathbf{s_k} \rangle , \quad\text{and}\quad\sigma_k = \langle  \mathbf{r_k}, \mathbf{s_k} \rangle.\]
Their equality in exact arithmetic holds, because, owing to the relation ${\mathbf{p_k} = \mathbf{r_k} + \beta_{k}\mathbf{p_{k-1}}}$ from Pipe-PR-CG, we can rewrite $\mu_k$ as 
\[\mu_k = \langle  \mathbf{p_k}, \mathbf{s_k} \rangle = \langle  \mathbf{r_k} + \beta_{k}\mathbf{p_{k-1}}, \mathbf{s_k} \rangle,\]
where the inner product $\langle  \mathbf{p_{k-1}}, \mathbf{s_k} \rangle$ is equal to zero. This holds because in exact arithmetic we have that $\mathbf{s_k} = \mathbf{Ap_k}$ and vectors $\mathbf{p_i}$, and $\mathbf{p_j}$ are $A$-orthogonal for $i \ne j$.
\par
Knowing this, we can define the $\mu$-gap as $\Delta_{\mu'_k} \equiv|\mu_k - \sigma_k|$, and  derive a bound for it.
We use the finite arithmetic relations 
\begin{equation}\label{420}
	 \mathbf{p_k} = \mathbf{r_{k}} + \beta_{k}\mathbf{p_{k-1}} + \delta_{\mathbf{p_k}}, 
	\quad \mu_k = \langle \mathbf{p_k}, \mathbf{s_k} \rangle + \delta_{\mu_k}, \quad \sigma_k = \langle \mathbf{r_k}, \mathbf{s_k} \rangle +\delta_{\sigma_k} .
\end{equation}
 Using these expressions, we can write 
\begin{align*}
	\mu_k - \sigma_k &= \langle  \mathbf{p_k}, \mathbf{s_k} \rangle - \langle  \mathbf{r_k}, \mathbf{s_k} \rangle + \delta_{\mu_k} - \delta_{\sigma_k}\\[3pt]
	&= \langle  \mathbf{r_k} + \beta_k \mathbf{p_{k-1}} + \delta_{\mathbf{p_k}} , \mathbf{s_k} \rangle - \langle  \mathbf{r_k} , \mathbf{s_k} \rangle + \delta_{\mu_k} - \delta_{\sigma_k}\\[3pt]
	&= \beta_k \langle  \mathbf{p_{k-1}} , \mathbf{s_k} \rangle + \langle \delta_{\mathbf{p_k}} , \mathbf{s_k} \rangle + \delta_{\mu_k} - \delta_{\sigma_k},
\end{align*}
which, by taking the norm of both sides, yields
\begin{equation}\label{421}
	\Delta_{\mu'_k} \le |\beta_k|\, |\langle  \mathbf{p_{k-1}} , \mathbf{s_k} \rangle| + ||\delta_{\mathbf{p_k}}|| \, ||\mathbf{s_k}|| + |\delta_{\mu_k}| + |\delta_{\sigma_k}|.
\end{equation}
Using \eqref{42_2} and \eqref{420}, $|\delta_{\mu_k}|$ and $|\delta_{\sigma_k}|$ can be bounded as
\begin{align}
	|\delta_{\mu_k}| &\le \epsilon \, n \,(||\mathbf{p_k}||\,||\mathbf{s_k}||),\label{422_1}\\
	|\delta_{\sigma_k}| &\le \epsilon \, n \,(||\mathbf{r_k}||\,||\mathbf{s_k}||), \label{422_2}
\end{align}
and, similarly, using \eqref{42_1} and \eqref{420}, we have
\begin{equation}\label{424}
	||\delta_{\mathbf{p_k}}|| \le \epsilon \,(||\mathbf{r_k}|| + 2|\beta_k|\,||\mathbf{p_{k-1}}||).
\end{equation}
\par
Unfortunately, rewriting the term $\langle  \mathbf{p_{k-1}} , \mathbf{s_k} \rangle$ in a way which involves norms of $\mathbf{p_{k-1}}$, $\mathbf{s_{k}}$, or some other variable, would greatly diminish tightness of the bound. Thus, by keeping this term and substituting \eqref{422_1}, \eqref{422_2}, and \eqref{424} into \eqref{421}, and dropping terms of order $O(\epsilon^2)$, we obtain that
\begin{equation*}
	\Delta_{\mu'_k} \lesssim  |\beta_k|\, |\langle  \mathbf{p_{k-1}} , \mathbf{s_k} \rangle| + \epsilon \, ||\mathbf{s_k}||\big(||\mathbf{r_k}|| + 2|\beta_k|\,||\mathbf{p_{k-1}}|| + n \, (||\mathbf{p_k}|| + ||\mathbf{r_k}||)\big).
\end{equation*}
{ Once again, by altering the derivation to take into account the preconditioning steps and quantities we can obtain a $\mu$-gap bound for the preconditioned version of Pipe-PR-CG:
\[
	\Delta_{\mu'_k} 
	 \le |\beta_k|\, |\langle  \mathbf{p_{k-1}} , \mathbf{s_k} \rangle| + \epsilon \, ||\mathbf{s_k}||\big(||\trk|| + 2|\beta_k|\,||\mathbf{p_{k-1}}|| + n \, (||\mathbf{p_k}|| + ||\trk||)\big).
\]
For these expressions, { three additional inner products are needed in each iteration to compute $\langle  \mathbf{p_{k-1}} , \mathbf{s_k} \rangle$, $\|\sk\|$, and $||\mathbf{p_k}||$}. There is no need to compute the norm $||\mathbf{p_{k-1}}||$ as we can simply keep its value from the previous iteration (or initialization). For future usage, let $B_{\mu'_k}$ denote the derived bound.}
\par
\begin{figure}[h]
	\makebox[\textwidth][c]{\includegraphics[clip, trim=0cm 0cm 0cm 1cm, width=.7\textwidth]{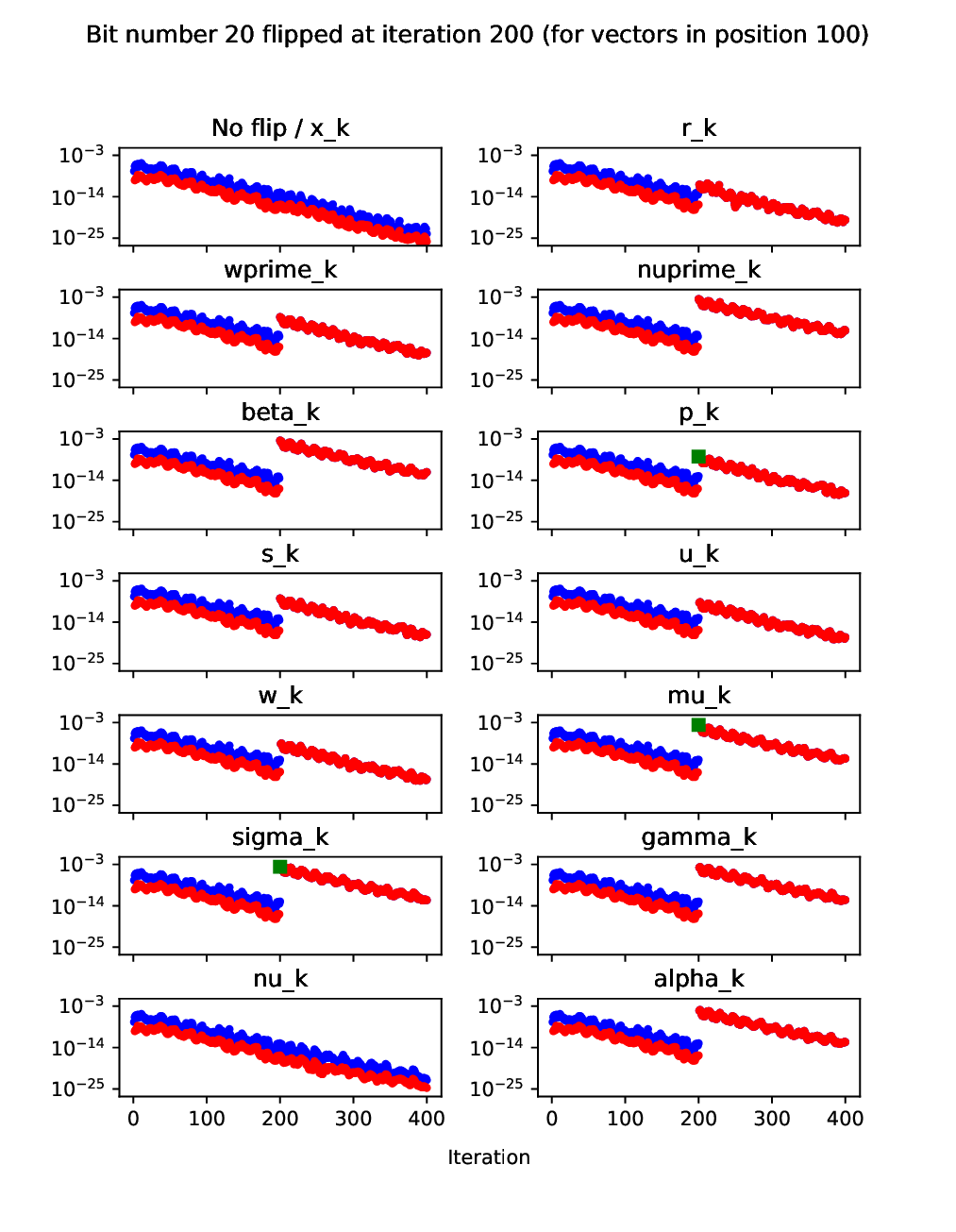}}
	\caption{$\mu$-gap (red) and $\mu$-gap bound (blue) graph, matrix \textit{bcsstm07}}
	\label{Fig7}
\end{figure}

Figure \ref{Fig7} depicts the $\mu$-gap and the $\mu$-gap bound using the same experimental setup as in the previous subsections. As before, if the bound is violated the gap marker turns into a square. However, this time the square is green in order to be easily visible, since for the $\mu$-gap it is always next to a cluster of other gap values.
\par
The first thing we observe is that the $\mu$-gap and bound are sensitive to flips in almost all variables. 
Another interesting fact is that the $\mu$-gap bound for many variables seemingly vanishes after the flip. However, this is caused by it being very close to the $\mu$-gap. It is also peculiar that neither the gap nor the bound return to their original level, but instead the values are permanently affected by the flip. This was not the case for the $\nu$- and the $\mathbf{w}$-gaps and bounds. The reason for this is the inner product $|\langle  \mathbf{p_{k-1}} , \mathbf{s_k} \rangle|$ present in the expression of the $\mu$-gap bound.
\par
It seems that the $\mu$-gap/bound method is able to detect flips in $\mathbf{p_k}$, $\mu_k$, and $\sigma_k$. For all three above-mentioned variables, the bound is always violated at the flip iteration. It is also interesting that the reason for this is that for these variables the $\mu$-gap ``jumps'' one iteration earlier than the $\mu$-gap bound. For all other variables this happens concurrently, either at the flip iteration or at the very next one.

There are still some variables for which we do not possess a detection method, namely $\mathbf{x_k}$, $\beta_k$, $\mathbf{s_k}$, and $\alpha_k$. However, a straightforward comparison of the gap values and the bound values is not the only way the detection can be done. As was mentioned before, and as can be seen by investigating Figure \ref{Fig7}, the $\mu$-gap and the $\mu$-gap bound are influenced by flips in almost all variables. On top of that, their values after the flip become very close. Therefore, we will to construct a detection method based on the difference of the $\mu$-gap and the $\mu$-gap bound. However, it somehow surprisingly turns out that their absolute difference, $|B_{\mu'_k} - \Delta_{\mu'_k}|$, steadily decreases despite flips (with the exception of a single jump for $\mathbf{p_k}$, $\mu_k$, and $\sigma_k$). Thus, we will employ the relative difference, $|B_{\mu'_k} - \Delta_{\mu'_k}|/B_{\mu'_k}$, instead.
\par
There are two reasons for ``normalizing'' the difference by the $\mu$-gap bound and not by the $\mu$-gap. The first one is that the $\mu$-gap can sometimes be zero. The second is that when no flips occur, the $\mu$-gap bound is guaranteed to be larger than the $\mu$-gap. Therefore, the ratio will be better ``normalized''.
\par
\begin{figure}[h]
	\makebox[\textwidth][c]{\includegraphics[clip, trim=0cm 0cm 0cm 1cm, width=.7\textwidth]{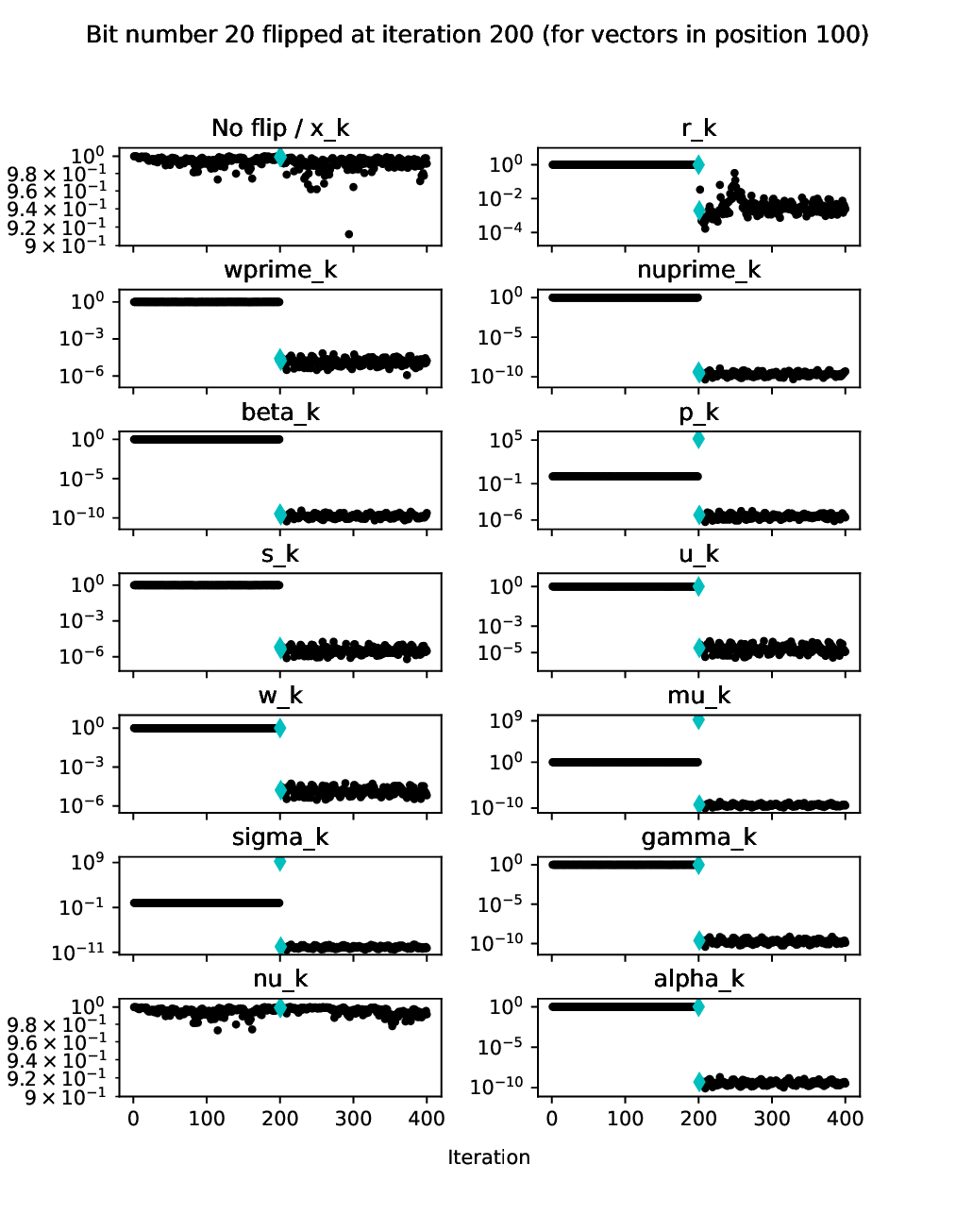}}
	\caption{Relative $\mu$-gap/bound difference, $|B_{\mu'_k} - \Delta_{\mu'_k}|/B_{\mu'_k}$, matrix \textit{bcsstm07}}
	\label{Fig8}
\end{figure}

Figure \ref{Fig8} shows the relative $\mu$-gap/bound difference. The cyan diamond markers show values in the flip iteration and in the one iteration after, i.e., when we would like to be able to detect the flip, so that we can roll back to a close previous state when the variables were still unaffected. The setup and data depicted are the same as in the previous figures of this section.

At first glance, we immediately observe that the effect of the flip is quite significant for all variables with the exception of $\nu_k$. Whether it is at the iteration of the flip, one later, or both, there is a significant jump of the value. In addition, for the three variables ($\mathbf{p_k}$, $\mu_k$, $\sigma_k$), where we were able to detect flips just by the bound violation, the $\mu$-gap/bound relative difference is greater than 1 at the flip iteration, a clear indication that a flip has occurred. 

It is also important to investigate at what level the values of the studied ratio $|B_{\mu'_k} - \Delta_{\mu'_k}|/B_{\mu'_k}$ are when no flips occur, because we need to set some threshold to determine whether to raise an alarm that a silent error has likely appeared or not. If we set the threshold too close to 1 we might incur many false positives. On the other hand, setting it too low may result in many false negatives. For some matrices, e.g., \textit{bcsstm07}, the values lie very close to 1. Unfortunately, for matrices with higher condition numbers the level of the ratio $|B_{\mu'_k} - \Delta_{\mu'_k}|/B_{\mu'_k}$ may not be so close to one as in Figure \ref{Fig8}. This is visible in Figure \ref{FigB4} in Appendix \ref{AppendixFigures} where the values of the ratio for matrix \textit{nos7} decrease to the level of $1\mathrm{e}{-4}$. The values there are altered by the flip, but it can be observed that for vector variables the diamond markers still largely remain in the same value range as when no bits were flipped. This indicates that choosing a suitable value of the threshold could be a rather complex and data dependent problem.

 However, there is only one vector variable, $\mathbf{s_k}$, which is not covered by the previous detection methods, and the relative $\mu$-gap/bound difference criterion seems to mostly work for this quantity. This can be also seen in Figure \ref{FigB4}. Nonetheless, we have also encountered cases where the reliability of this criterion for $\mathbf{s_k}$ is borderline.

In summary, the greatest strength of the relative $\mu$-gap/bound difference approach is that it encompasses almost all of the Pipe-PR-CG variables, albeit with some above-mentioned data-related exceptions. However, a disadvantage is that, unlike in the case of the bound violation methods, there is nothing to directly compare the values to, so we have to set a detection threshold.

\subsection{Summary of detection methods}
We have investigated four methods for silent error detection in Pipe-PR-CG. We summarize their efficacy in terms of whether or not they have the potential to reliably detect bit flips in a given variable in Table \ref{table1}, with rows representing each of the Pipe-PR-CG variables. 
The symbol \checkmark denotes that the method is able to reliably detect flips in the variable, while $\circ$ denotes that the method is, for the given variable, somehow functional, but either not in all cases or there are some specific circumstances under which it is not able to detect the injected error, e.g., the $\nu$-gap/bound criterion not working for sign flips in the residual vector $\mathbf{r_k}$. 
\par
We can see from the table that the only variable which is not covered by any of the detection methods is the solution vector $\mathbf{x_k}$. The reason for this is that $\mathbf{x_k}$ appears only in its own relation, and thus it does not affect any other variable. Therefore, for silent error detention in $\mathbf{x_k}$, a redundancy approach is unfortunately needed. Besides that, we may not posses a robust detection method for flips in $\mathbf{s_k}$. Nonetheless, in the worst case, the redundancy approach can be applied here as well if it turns out that detection by the relative $\mu$-gap/bound difference is truly unreliable. For all other variables we should be able to detect silent errors reasonably well.

{We note that our detection criteria require the computation of additional inner products. These, however, can still be computed at the same time (in a single synchronization point) as the other inner products in Pipe-PR-CG. In addition, to reliably compute the detection criteria, these additional inner products should be computed redundantly, since errors could also occur during their computation. This can still be accomplished in the same synchronization point, and thus does not diminish the potential for pipelining. We note that the bounds themselves should also be computed redundantly; these are simple computations with scalar quantities and thus this should not be a significant cost.  }

We note that if the computation of additional inner products and constants necessary for the bounds would in some case be too expensive, it could also possible to construct a set of detection criteria based just on the values of the gaps alone. For instance, we could monitor a moving average of gap value differences between iterations. Nonetheless, it is important to keep in mind that for this approach to function we must somehow deal with iterations where any of the gaps are zero. It would also require us to set some threshold to determine when to raise an alarm that a silent error has likely occurred.

\begin{table}[h]\
	\centering
    \footnotesize
	\renewcommand{\arraystretch}{1.2}
	\begin{tabular}{ |c||c|c|c|c| } 
		\hline
		Variable& $\nu$-gap/bound & $\mathbf{w}$-gap/bound &  $\mu$-gap/bound&$|B_{\mu'_k} - \Delta_{\mu'_k}|/B_{\mu'_k}$\\
		\hline
		$\mathbf{x_k}$ & & & &\\
		\hline
		$\mathbf{r_k}$ & $\circ$ & \checkmark & &$\circ$\\ 
		\hline
		$\mathbf{w'_k}$ & & \checkmark & &$\circ$\\ 
		\hline
		$\nu'_k$ & \checkmark & & &\checkmark\\ 
		\hline
		$\beta_k$ & & & &\checkmark\\ 
		\hline
		$\mathbf{p_k}$ & & & \checkmark &\checkmark\\ 
		\hline
		$\mathbf{s_k}$ & & & &$\circ$\\ 
		\hline
		$\mathbf{u_k}$ & & \checkmark & &$\circ$\\ 
		\hline
		$\mathbf{w_k}$ & & \checkmark & &$\circ$\\ 
		\hline
		$\mu_k$ & & & \checkmark &\checkmark\\ 
		\hline
		$\sigma_k$ & \checkmark & & \checkmark &\checkmark\\ 
		\hline
		$\gamma_k$ & \checkmark & & &\checkmark\\ 
		\hline
		$\nu_k$ & \checkmark & & &\\ 
		\hline
		$\alpha_k$ & & & &\checkmark\\ 
		\hline
	\end{tabular}
	\renewcommand{\arraystretch}{1}
	\caption{Efficacy of detection methods for Pipe-PR-CG variables} 
 \label{table1}
\end{table}

{
\subsection{Effects of preconditioning}
Throughout this section, we have briefly mentioned what the $\nu$-gap bound, the $\mathbf{w}$-gap bound, and the $\mu$-gap bound look like for the Pipe-PR-CG algorithm with preconditioning. Let us now elaborate more thoroughly about what happens when  preconditioning is included. First, as was already stated, the bounds change to the following: 
\begin{align}
    \Delta_{\nu_k'} &\lesssim \left( \frac{21+6n}{2} \right) \left(\|\rrkm\|^2 + \|\trkm\|^2 + \|\rrk\|^2 + \|\trk\|^2 \right),\\
    \Delta_{w'_k} &\lesssim  2(c+3)\epsilon||\mathbf{A}||\, (||\mathbf{\tilde{r}_{k-1}}||+||\mathbf{\tilde{r}_{k}}||\label{wboundpred}\\
    \Delta_{\mu'_k} 
	 &\lesssim |\beta_k|\, |\langle  \mathbf{p_{k-1}} , \mathbf{s_k} \rangle| + \epsilon \, ||\mathbf{s_k}||\big(||\trk|| + 2|\beta_k|\,||\mathbf{p_{k-1}}|| + n \, (||\mathbf{p_k}|| + ||\trk||)\big). 
\end{align}

The set of detection criteria based on violation of these three bounds and the ratio $|B_{\mu'_k} - \Delta_{\mu'_k}|/B_{\mu'_k}$ was shown to encompass all Pipe-PR-CG variables appearing in the unpreconditioned case. However, in the preconditioned case, there are several additional variables, namely, $\trk$, $\mathbf{\tilde{w}'_k}$, $\tsk$, $\mathbf{\tilde{u}_k}$, and $\mathbf{\tilde{w}_k}$.

In our experience, not all of these ``preconditioning-exclusive'' variables cause jumps in the $\nu$-gap, the $\mathbf{w}-$gap, or the $\mu$-gap. Therefore, one must introduce additional criteria for the preconditioned case. For instance, one can derive a bound for the so-called $\mathbf{\tilde{w}_k}$-gap, $\Delta_\twkp \equiv \|\twk - \twkp\|$. The analysis is generally straightforward, although here we will need to capture the action of the preconditioner on a vector. We assume, for example, we have the relation 
\begin{equation*}
      \twk = \M\wk + \delta_\twk.
\end{equation*}
Clearly, a bound on $\delta_\twk$ will depend on the particular form of the preconditioner used and how it is applied, e.g., whether it involves triangular solves as with a Cholesky factor, or whether it can be applied via a matrix-vector product as with a sparse approximate inverse. Fortunately, these two mentioned computations are widely applicable to a range of commonly-used algebraic preconditioners, so the additional analysis that must be done for many practical cases is minimal. We keep our analysis general here and assume only that the application of $\M$ can be written in a form such that 
\[
\delta_\twk = \delta_\M \wk, \quad\text{where}\quad |\delta_\M|\leq f(n) \epsilon E
\]
where $E$ has positive entries and $f(n)$ is polynomial in $n$, and thus we will bound, e.g., 
\begin{align}
\|\delta_\twk \| &\leq f(n)\epsilon\|\delta_\M\|\|\wk\| \nonumber\\
&\leq f(n)\epsilon\|\delta_\M\|\|A\|\|\trk\|+O(\epsilon^2) \label{deltatwk}.
\end{align}
Using this result and following standard rounding error analysis, we can obtain the $\mathbf{\tilde{w}_k}$-gap bound 
\begin{equation}
\Delta_\twkp \lesssim \Big[2(c+3) \|\M\| +2f(n) \|\delta_\M\|\Big] \epsilon \|\A\|\big(\|\trkm\|+\|\trk\| \big). \label{wtildeboundpred}
\end{equation}
Note that this bound requires, in addition to an estimate of $\|\A\|$, an estimate of the norm of the inverse of the preconditioner, as well as an estimate of $\|\delta_\M\|$, which is dependent on the way in which $\M$ is applied. In our experiments, we have found that using the measured norms of $\M$ and $\delta_\M$ in the bound leads to a very large overestimate of the error; this means that even if there is a significant jump in the value of $\|\twk-\twkp\|$, it may still not violate the bound. In many practical scenarios, we believe it will suffice to use the bound for the $\mathbf{w}$-gap in (\ref{wboundpred}) in place of (\ref{wtildeboundpred}) in detecting errors (essentially setting the norm of $\M$ to 1 and $\delta_\M$ to $0$). The only risk is that this could result in too many false positives.

Based on our initial experimentation, we believe that the detection of silent errors for Pipe-PR-CG algorithm with preconditioning is a very challenging problem, which is highly dependent on the data and choice of particular preconditioner. Because of this, we leave it as an important topic for future research.
}

\section{Fault-tolerant Pipe-PR-CG}\label{sec:ftpipeprcg}
We now aim to test how well our criteria can detect silent errors. We examine the performance of the criteria on a large sample of test runs, both with and without bit flips. The experiment was performed for each of the Pipe-PR-CG variables with the exception of $\mathbf{x_k}$, since none of the methods work for this variable. It is also worth noting that the experiment was performed for each of the variables separately, so that eventual outliers can be identified more easily. 

For all variables, the testing was done using the first eight matrices in Table \ref{tab:matrices}. 
For each variable and each of these matrices, 800 runs with a single bit flip and 200 without a bit flip were performed. The bit number was chosen randomly from 1 to 64. For vectors, the flip occurred in a random index from $[1,n]$, where $n$ is the problem dimension. The flip iteration $\tau$ was chosen randomly from $0.1\varphi$ to $0.9\varphi$, where $\varphi$ is the number of iterations needed to converge for the given matrix and right-hand side when no flips occur. This was computed before each tainted run. The stopping criterion was always, for both untainted as well as tainted runs, such that it must hold that $||\mathbf{r_k}||/||\mathbf{b}|| < 1\mathrm{e}{-10}$. The right-hand side $\mathbf{b}$ was a random vector from a uniform distribution over $[0,1)$. A run tainted by a bit flip was deemed as convergent if it reached the stopping criteria within $1.5\varphi$ iterations. The initial guess $\mathbf{x_0}$ was always a vector of all zeros. Runs with an overflow error were not counted, because such errors are no longer silent. For this reason, the total number of runs recorded for each matrix is slightly lower than the above-stated 1000 performed. When possible, the norms appearing in the code were calculated using the already computed quantities, such as $\nu_k$ for $||\mathbf{r_k}||$ or $\gamma_k$ for $||\mathbf{s_k}||$. The only exception was the norm $||\mathbf{r_k}||$ used for the stopping criterion which was not computed utilizing $\nu_k$, so that the convergence can be evaluated independently from the detection criteria.
\par
Let us now recall what our four detection methods are.  If any of the inequalities
\begin{itemize}
	\item $|\nu_k - \nu_k'| > \epsilon \, (21+6n)(||\mathbf{r_{k-1}}||^2 + ||\mathbf{r_{k}}||^2)$,
	\item $||\mathbf{w_k - w'_k}|| > 2(c+3)\epsilon  ||\mathbf{A}||\, \big( ||\mathbf{r_k}|| + ||\mathbf{r_{k-1}}|| \big)$,
	\item $\Delta_{\mu'_k} \equiv |\mu_k - \sigma_k| \newline > |\beta_k|\, |\langle  \mathbf{p_{k-1}} , \mathbf{s_k} \rangle| + \epsilon \, ||\mathbf{s_k}||\big(||\mathbf{r_k}|| + 2|\beta_k|\,||\mathbf{p_{k-1}}|| + n \, (||\mathbf{p_k}|| + ||\mathbf{r_k}||)\big) \eqqcolon B_{\mu'_k}$,
	\item $|B_{\mu'_k} - \Delta_{\mu'_k}|/B_{\mu'_k} < \text{threshold}$,
\end{itemize}
held, an alarm was raised. We call these four criteria working together a \textit{detection set}. In the experiment, there were two detection sets with a different threshold for the relative $\mu$-gap/bound difference criterion to evaluate how the detection behavior changes with the threshold. The threshold values were chosen to be $5\mathrm{e}{-1}$ and $1\mathrm{e}{-4}$, as these values seemed to be close to the lower limit of values of the relative $\mu$-gap/bound difference, $|B_{\mu'_k} - \Delta_{\mu'_k}|/B_{\mu'_k}$, for matrices \textit{1138\_bus} and \textit{nos7} when no flips occur. The two detection sets with different thresholds for the fourth relative $\mu$-gap/bound difference-based criterion were both evaluated simultaneously during each run. Therefore, we can directly compare how they perform for identical data.
\par
The sequence of steps in the experiment was following:
\begin{enumerate}
	\item For untainted runs, a right-hand side vector $\mathbf{b}$ is generated, then the computation is performed and it is noted whether an alarm was raised. The two detection sets using the two different thresholds for the $|B_{\mu'_k} - \Delta_{\mu'_k}|/B_{\mu'_k}$ criterion are both checked during the run and they each posses their own alarm.
	\item For tainted runs, a right-hand side vector $\mathbf{b}$ is generated, and subsequently an untainted computation is performed to obtain the number of iterations $\varphi$ needed to converge. Afterwards, the flip iteration $\tau$, bit number, and vector flip index are generated. Then, a tainted run is performed with these inputs. As was the case for the untainted runs, both detection sets are monitored during the computation. Once again, this is done so that they can be better compared against each other, since they examine the same data. If during the run one of the two detection sets raises its alarm, it is noted at what iteration $\rho_i$, $i \in \{1,2\}$, that first was. Later alarms are not taken into account. Besides the first alarm iterations, it is also monitored whether the run converged within the $1.5\varphi$ iterations or not.
\end{enumerate}

Runs were sorted into one of six categories:
\begin{itemize}
	\item true positive (tp): Flip occurred, did prevent convergence, alarm was raised.
	\item special positive (sp): Bit flip occurred, did not prevent convergence, alarm was raised.
	\item false positive (fp): No flip occurred, alarm was raised.
	\item true negative (tn): No flip occurred, no alarm was raised.
	\item special negative (sn): Flip occurred, did not prevent convergence, no alarm raised.
	\item false negative (fn): Flip occurred, did prevent convergence, no alarm raised.
\end{itemize}
The same categorization was used in \cite{MeurantCG}. As mentioned above, if in a tainted run any of the four criteria within the detection sets raised an alarm it was noted at what iteration ($\rho_1$ for the first detection set and $\rho_2$ for the second detection set) this first occurred. The values $\rho_i$ were initialized to $\infty$. Thus, if the detection set did not raise an alarm during the run its $\rho_i$ value remained $\infty$ after the computation had concluded. To classify runs, we compared the flip iteration $\tau$ and $\rho_i$, receiving the following options: 

\begin{enumerate}
	\item If $\rho_i < \tau$, we count this as false positive,
	\item If $\rho_i \in \{\tau, \tau + 1\}$, we count this as true/special positive based on whether the run converged,
	\item If $\rho_i > \tau + 1$, we count this as false/special negative based on whether the run converged.
\end{enumerate}
Runs without a bit flip were categorized either as true negative or false positive based on whether or not an alarm was raised.
\par
The output of the experiment is presented below in Table \ref{table2}, which contains the sum of all runs over the individual variables. 
The most crucial result is that we were generally able to detect an overwhelming majority of bit flips which would destroy convergence. 
Important also is the fact that there was no variable which would stick out as seriously problematic for our detection methods. Moreover, the number of false positive runs was (for both thresholds) very close for all variables.
\par
However, there were differences in the overall number of detected errors which would not destroy convergence (sp/sn). The one variable which stands out is $\mathbf{s_k}$, for which the number of special negative runs was considerably larger than for any other variable. This was most likely caused by the fact that for $\mathbf{s_k}$ only the relative $\mu$-gap/bound difference criterion works, and even then, it is not fully reliable, as previously mentioned. Nonetheless, a majority of the undetected flips were special negative, thus they did not destroy convergence, and it can be seen that the number of false negatives is for $\mathbf{s_k}$ quite acceptable. 
Notable also is the fact that in the case of scalar variables for which one of the gap/bound detection methods works ($\nu'_k$, $\mu_k$, $\sigma_k$, $\gamma_k$, and $\nu_k$) we were able to detect a large portion of the convergence-preserving silent errors. 
\par
For most matrices, there were no or almost no false positives. Notable outliers are \textit{nos7} and \textit{1138\_bus}. Interestingly, these two matrices along with \textit{bcsstk05} and \textit{bcsstm07} are all in the upper half of our sample when it comes to condition number. This leads to the conclusion that the threshold value should be ideally chosen proportionally to the condition number of the matrix. We note that false positive detections are caused only by the criterion utilizing the relative difference of the $\mu$-gap and the $\mu$-gap bound. The three bound violation criteria raise the alarm only when a bit flip truly occurs.

\begin{table}[h]\
	\footnotesize
	\centering
	\renewcommand{\arraystretch}{1.1}
	\begin{tabular}{c c c c c c c c } 
		matrix&threshold&tp& sp & fp & tn & sn & fn \\
		\hline
		\textit{1138\_bus}&$5\mathrm{e}{-1}$&1532 &	4125	&2681	&1652	&2844	&2 \\
		&$1\mathrm{e}{-4}$&1732	&4509&	0&	2600&	3989&	6 \\
		\hline
		\textit{bcsstm07}&$5\mathrm{e}{-1}$&939&	6342&	2&	2600&	2960&	3 \\
		&$1\mathrm{e}{-4}$&940	&5710	&0&	2600&	3593	&3 \\
		\hline
		\textit{bundle1}&$5\mathrm{e}{-1}$&833&	5213&	0	&2600&	4088&	1 \\
		&$1\mathrm{e}{-4}$&833&	4554&	0&	2600&	4747&	1 \\
		\hline
		\textit{wathen120}&$5\mathrm{e}{-1}$&771&	5331	&0	&2600	&4113&	2 \\
		&$1\mathrm{e}{-4}$&771&	4769&	0&	2600&	4675&	2 \\
		\hline
		\textit{bcsstk05}&$5\mathrm{e}{-1}$&2991&	4793&	348	&2324&	2314&	5 \\
		&$1\mathrm{e}{-4}$&2999&	4186&	0&	2600&	2980&	10 \\
		\hline
		\textit{gr\_30\_30}&$5\mathrm{e}{-1}$&964&	6479&	0	&2600	&2783&	2 \\
		&$1\mathrm{e}{-4}$&964	&5900&	0	&2600&	3362&	2\\
		\hline
		\textit{nos7}&$5\mathrm{e}{-1}$&0&	0&	12823&	0&	0&	0 \\
		&$1\mathrm{e}{-4}$&1381	&2881&	3843&	1229&	3485&	4 \\
		\hline
		\textit{crystm01}&$5\mathrm{e}{-1}$&799&	5896&	0	&2600	&3564&	1 \\
		&$1\mathrm{e}{-4}$&799&	5292&	0&	2600&	4168	&1 \\
		\hline
		\hline
		$\sum$&$5\mathrm{e}{-1}$&8829&       38179     &  15854     &  16976  &     22666 &         16\\
		&$1\mathrm{e}{-4}$&10419    &   37801  &      3843    &   19429& 30999     &     29\\
	\end{tabular}
	\renewcommand{\arraystretch}{1}
	\caption{Detection performance, sum over all variables} 
 \label{table2}
\end{table}

The combination of our detection methods is capable of reliably detecting the majority of silent errors which would destroy the convergence of Pipe-PR-CG. The problem at hand is now how to correct these errors. The approach we present here is to perform a so-called \textit{rollback} when the alarm is raised. Rolling back essentially means to ``return'' the computation to an uncorrupted state before the detected silent error has occurred. Therefore, in our case, to \textit{recover} the computation we have to ``return'' two iterations back, since our detection methods raise the alarm either at the iteration when the error has occurred or one iteration later. This recovery approach was previously proposed in other studies investigating the detection and correction of silent errors, e.g., in \cite{MeurantCG}.
\par

Algorithm \ref{FT-Pipe-PRCG} gives the resulting Fault-Tolerant Pipelined Predict-and-Recompute Conjugate Gradient algorithm (FT-Pipe-PR-CG), which includes the detection and correction of silent errors. 
We have also explicitly incorporated here the redundancy detection approach for $\mathbf{x_k}$. Inputs of the algorithm include, aside from the standard problem data $\mathbf{A}$, $\mathbf{b}$, and $\mathbf{x_0}$, the constants $||\mathbf{A}||$, $n$, $\epsilon$, and $c$, necessary for evaluation of the bounds, and the threshold value $T$ used in the $|B_{\mu'_k} - \Delta_{\mu'_k}|/B_{\mu'_k}$ criterion.
\par
If the alarm is raised we call a Recover() procedure (Algorithm \ref{Recover}) to perform the two iteration rollback and we mark the iteration $k+2$ as ``corrected''; otherwise, if the detection was falsely positive, the alarm would be raised again indefinitely every two iterations. Although it may theoretically happen that another silent fault appears in an iteration marked as corrected, the probability of this is rather low, since silent errors are a rare event \cite{ElliottIEEE}. However, in cases where there is an extremely large number of false positive detections, this may be a problem. We later propose a remedy.
\par 
In FT-Pipe-PR-CG it is necessary to compute three additional inner products, $\langle \mathbf{w_k - w'_k}, \mathbf{w_k - w'_k} \rangle$, $\langle \mathbf{p_{k-1}}, \mathbf{s_k} \rangle$, and $\langle \mathbf{p_k}, \mathbf{p_k} \rangle$, which do not appear in the basic Pipe-PR-CG algorithm, for the sake of our detection criteria. 
{ Aside from $\langle \mathbf{w_k - w'_k}, \mathbf{w_k - w'_k} \rangle$, it is possible to couple their computation with the other inner products in line \ref{FTPPRCG-ips} of Algorithm \ref{FT-Pipe-PRCG}. The problem with the $\mathbf{w}$-gap is that the variable $\mathbf{w_k}$ necessary for its computation is communicated among the computational units during the singular global synchronization point. However, computation of the inner product $\langle \mathbf{w_k - w'_k}, \mathbf{w_k - w'_k} \rangle$ requires a global communication as well. Luckily, there is a solution that will allow us to retain only one global synchronization point. We can perform the necessary communication for the inner product as well as the evaluation of the $\mathbf{w}$-gap/bound criterion in the next iteration. Thus, we will be partly performing the silent error detection for iteration $k$ during iteration $k+1$. This means that even though our criteria are able to detect silent errors based on variables from at worst one iteration after the error occurs, an error at iteration $k$ may be corrected two iterations later ($k+2$), albeit based on data from iteration $k+1$. Note that it is also possible to evaluate all the other detection criteria along with the $\mathbf{w}$-gap/bound ``one iteration later''. However, in our opinion it is better to evaluate them as soon as possible, so that in case an error is detected, we do not have to redundantly compute an extra ``faulty'' iteration.}

The advantage of the rollback correction approach is that the algorithm is able to universally recover from any detected silent error, no matter what variable it occurred in. The disadvantage is that we need to allocate extra memory for storing the variables from iterations $k-2$, $k-3$, { and $k-4$}.

{
\begin{algorithm}[h]
	\setstretch{0.90}
	\caption{Recover procedure of FT-Pipe-PR-CG}
	\label{Recover}
	\begin{algorithmic}[1]
		\Procedure {Recover}{$\cdot$}
		\State $\mathbf{x_k} = \mathbf{x_{k-3}}, \quad \mathbf{x_{k-1}} = \mathbf{x_{k-4}}$, $\quad \mathbf{r_k} = \mathbf{r_{k-3}}, \, \, \quad \mathbf{r_{k-1}} = \mathbf{r_{k-4}}$
		\State $\mathbf{w'_k} = \mathbf{w'_{k-3}}, \: \, \, \mathbf{w'_{k-1}} = \mathbf{w'_{k-4}}$, $\quad \nu'_k = \nu'_{k-3}, \quad \, \; \nu'_{k-1} = \nu'_{k-4}$
		\State $\beta_k = \beta_{k-3}, \quad \; \beta_{k-1} = \beta_{k-4}$, $\quad \mathbf{p_k} = \mathbf{p_{k-3}}, \; \: \, \; \mathbf{p_{k-1}} = \mathbf{p_{k-4}}$
		\State $\mathbf{s_k} =\mathbf{s_{k-3}}, \, \, \quad \mathbf{s_{k-1}} =\mathbf{s_{k-4}}$, $\quad \mathbf{u_k} = \mathbf{u_{k-3}}, \; \, \, \; \mathbf{u_{k-1}} = \mathbf{u_{k-4}}$
		\State $\mathbf{w_k} = \mathbf{w_{k-3}}, \; \; \mathbf{w_{k-1}} = \mathbf{w_{k-4}}$, $\quad \mu_k = \mu_{k-3}, \quad \: \mu_{k-1} = \mu_{k-4}$ 
		\State $\sigma_k = \sigma_{k-3}, \quad \: \sigma_{k-1} = \sigma_{k-4}$, $\quad \gamma_k = \gamma_{k-3}, \quad \: \gamma_{k-1} = \gamma_{k-4}$ 
		\State $\nu_k = \nu_{k-3}, \quad \, \: \nu_{k-1} = \nu_{k-4}$, $\quad\alpha_k = \alpha_{k-3}, \quad  \alpha_{k-1} = \alpha_{k-4}$
		\EndProcedure
	\end{algorithmic}
\end{algorithm}

\begin{algorithm}[h]
	\setstretch{1.0}
	\caption{Fault-Tolerant Pipelined Predict-and-Recompute Conjugate \newline Gradient: FT-Pipe-PR-CG}
	\label{FT-Pipe-PRCG}
	\begin{algorithmic}[1]
		\Procedure {FT-Pipe-PR-CG}{$\mathbf{A}, \mathbf{b}, \mathbf{x_0}$, $||\mathbf{A}||$, $n$, $\epsilon$, $c$, $T$, $tol$}
		\State INITIALIZE()
		\While {$ \|\mathbf{r_k}\|/\|\mathbf{b}\| > tol$}
		\State $\mathbf{x_k} = \mathbf{x_{k-1}} + \alpha_{k-1}\mathbf{p_{k-1}},\quad\mathbf{\hat{x}_k} = \mathbf{x_{k-1}} + \alpha_{k-1}\mathbf{p_{k-1}}$ 
        \State \textbf{if} $\mathbf{x_k \ne \hat{x}_k}$ \textbf{then go to 4} 

		\State $\mathbf{r_k} = \mathbf{r_{k-1}} - \alpha_{k-1}\mathbf{s_{k-1}},\quad\mathbf{w'_k} = \mathbf{w_{k-1}} - \alpha_{k-1}\mathbf{u_{k-1}}$
		\State $\nu'_k = \nu_{k-1} - 2\alpha_{k-1} \sigma_{k-1}+\alpha^2_{k-1}\gamma_{k-1},\quad\beta_k = \nu'_k/\nu_{k-1}$
		\State $\mathbf{p_k} = \mathbf{r_k} + \beta_{k}\mathbf{p_{k-1}},\quad\mathbf{s_k} = \mathbf{w'_k}+\beta_{k}\mathbf{s_{k-1}},\quad\mathbf{u_k} = \mathbf{A s_k},\quad\mathbf{w_k} = \mathbf{A r_k}$
		\State $\mu_k = \langle \mathbf{p_k}, \mathbf{s_k} \rangle,\quad \sigma_k = \langle \mathbf{r_k}, \mathbf{s_k} \rangle,\quad \gamma_k = \langle \mathbf{s_k}, \mathbf{s_k} \rangle,\quad \nu_k = \langle \mathbf{r_k}, \mathbf{r_k} \rangle$ \label{FTPPRCG-ips}
		\State $\alpha_k = \nu_k/\mu_k$

		\State $\Delta_{\nu'_k} = |\nu_k - \nu'_k|,\quad{\Delta_{\mathbf{w}'_{k-1}} = ||\mathbf{w_{k-1} - w'_{k-1}}||},\quad\Delta_{\mu'_k} = |\mu_k - \sigma_k|$
		\State $B_{\nu'_k} = \epsilon \, (21+6n)(||\mathbf{r_{k-1}}||^2 + ||\mathbf{r_{k}}||^2)$
		\State $B_{\mathbf{w}'_k} = \epsilon \, ||\mathbf{A}||\, ( (c+3)||\mathbf{r_k}|| + (c+4)||\mathbf{r_{k-1}}|| + (c+2)\,|\alpha_{k-1}|\,||\mathbf{s_{k-1}}|| )$
		\State $B_{\mu'_k} = |\beta_k|\, |\langle  \mathbf{p_{k-1}} , \mathbf{s_k} \rangle| + \epsilon \, ||\mathbf{s_k}||(||\mathbf{r_k}|| + 2|\beta_k|\,||\mathbf{p_{k-1}}|| + n \, (||\mathbf{p_k}|| + ||\mathbf{r_k}||))$
		\NoThenIf{$\Delta_{\nu'_k}>B_{\nu'_k}$ \textbf{or}  {$\Delta_{\mathbf{w}'_{k-1}}>B_{\mathbf{w}'_{k-1}}$} \textbf{or} $\Delta_{\mu'_k}>B_{\mu'_k}$ \textbf{or} $|B_{\mu'_k} - \Delta_{\mu'_k}|/B_{\mu'_k} < T$} 
		\NoThenIf{\textbf{$k$ is not marked as corrected}}
		\State \textbf{Recover()}, {\textbf{mark $k$ as corrected, $k=k-3$}}
		\EndIf
		\EndIf
		\EndWhile
		\EndProcedure
	\end{algorithmic}
\end{algorithm}
}

\section{An Adaptive Threshold Approach}\label{sec:adaptive}
FT-Pipe-PR-CG should be able to reliably detect and correct the majority of silent errors significant for convergence. However, in the previous section (e.g., Table \ref{table2}) we have observed that for some matrices there were many runs which resulted in false positives. Moreover, the experiment was categorizing the runs based only on the \textit{first} raising of the alarm. Hence, in some problematic cases there may potentially be a large number of false positive detections during the computation. This would cause us to perform many extra iterations due to the rollback recovery.
However, we have also seen that the number of false positives decreases with the value of the threshold parameter. 
\par
We will now instead \emph{adapt} the value of the threshold $T$ during the run of the algorithm to reflect how many times the alarm was raised. As was mentioned, silent errors are rather rare events, so if the alarm is raised many times we can safely assume that in most cases we did not truly detect a fault. In such a situation it may be beneficial to lower the value of the threshold $T$ to reduce the number of false positive detections by the relative $\mu$-gap/bound difference criterion. This idea is presented below in Algorithm \ref{AFT-Pipe-PRCG} as the adaptive fault-tolerant Pipe-PR-CG (AFT-Pipe-PR-CG). Here we multiply $T$ by an adaptation parameter $a \in (0,1)$ each time the alarm is raised by the $|B_{\mu'_k} - \Delta_{\mu'_k}|/B_{\mu'_k}$ criterion. Note that it is also possible to increase the threshold when there is a large number of iterations without any alarm; in this case, it is necessary to set some upper limit for $T$.
\par 
Not only does the adaptive algorithm potentially greatly reduce the number of false positive detections, it also allows us to eliminate the iteration marking. In FT-Pipe-PR-CG, if the alarm was raised at some iteration $k$ we have marked iteration $k+2$ as corrected, so that the procedure cannot get stuck in a loop. However, this can be caused only by the relative $\mu$-gap/bound difference criterion.
As was mentioned earlier, the three bound violation criteria raise the alarm only when a silent error truly occurs, i.e., they do not cause false positive detections. 
Thus, the procedure cannot get stuck because one of these criteria will indefinitely force a recovery in some iteration. The relative $\mu$-gap/bound difference criterion could do this, but now, each time this method raises the alarm the threshold is lowered. Therefore, eventually it will hold that $T<|B_{\mu'_k} - \Delta_{\mu'_k}|/B_{\mu'_k}$, and the procedure will continue.
\par
\begin{figure}[h]
	\makebox[\textwidth][c]{\includegraphics[clip, trim=0cm 0cm 0cm 1cm, width=.7\textwidth]{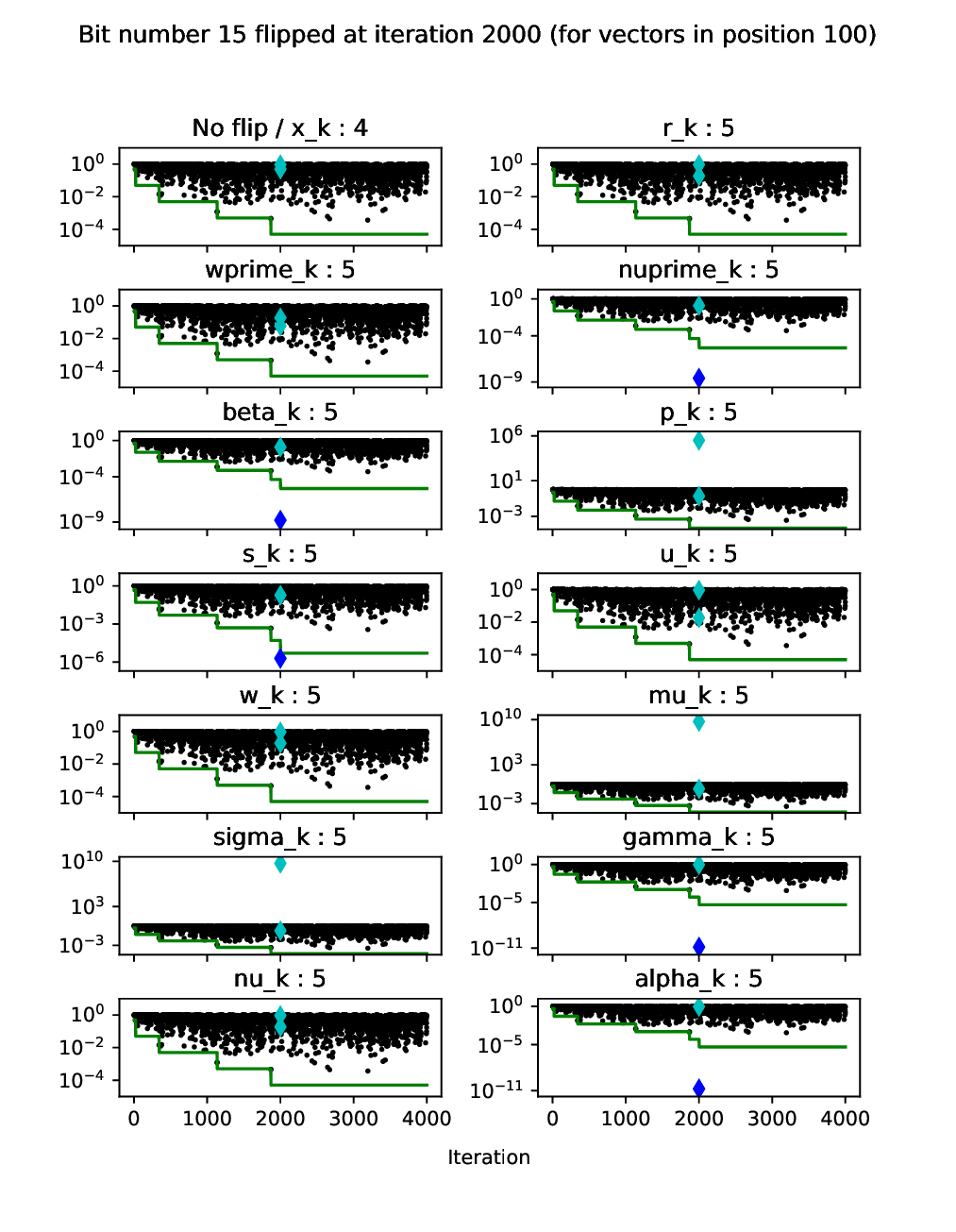}}
	\caption{Adaptive threshold $T$ for $a=0.1$ (green) and the relative $\mu$-gap/bound difference, $|B_{\mu'_k} - \Delta_{\mu'_k}|/B_{\mu'_k}$, (black), matrix \textit{nos7}}
	\label{Fig9}
\end{figure}

Figure \ref{Fig9} shows the process of threshold adaptation for the matrix \textit{nos7} and right-hand side $\mathbf{b = e}$. The adaptivity parameter $a$ was set to $0.1$. The initial threshold value was $5\mathrm{e}{-1}$, the higher value used in the detection performance experiment earlier in this section. Next to the variable names, it is noted how many recoveries, i.e., detections, there were in total during the computation. This number also includes alarms raised after the bit flip by criteria other than the threshold violation by the relative $\mu$-gap/bound difference. Therefore, for some variables, e.g., $\mathbf{w'_k}$, the number of threshold adaptations was one less than the number of total detections indicated in the figures. A violation of the threshold by the $|B_{\mu'_k} - \Delta_{\mu'_k}|/B_{\mu'_k}$ ratio in the flip iteration or one iteration later is denoted by the diamond marker turning dark blue.
\par
We observe that AFT-Pipe-PR-CG seems to be able to suitably adapt the threshold, so that the number of false positive detections is reduced, but at the same time the reliability of the detection is not destroyed. Notable also is that the ratio $|B_{\mu'_k} - \Delta_{\mu'_k}|/B_{\mu'_k}$ no longer ``jumps down'' as we have observed in Figure \ref{Fig8}. This is because of the recovery procedure.

{
\begin{algorithm}[h]
	\setstretch{1.20}
	\caption{Adaptive Fault-Tolerant Pipelined Predict-and-Recompute \newline Conjugate Gradient: AFT-Pipe-PR-CG}
	\label{AFT-Pipe-PRCG}
	\begin{algorithmic}[1]
		\Procedure {AFT-Pipe-PR-CG}{$\mathbf{A}, \mathbf{b}, \mathbf{x_0}$, $||\mathbf{A}||$, $n$, $\epsilon$, $c$, $T$, $a$, $tol$}
		\State INITIALIZE()
		\While {$ \|\mathbf{r_k}\|/\|\mathbf{b}\| > tol$}
		\State $\mathbf{x_k} = \mathbf{x_{k-1}} + \alpha_{k-1}\mathbf{p_{k-1}}, \quad \mathbf{\hat{x}_k} = \mathbf{x_{k-1}} + \alpha_{k-1}\mathbf{p_{k-1}}$ \State{\textbf{if} $\mathbf{x_k \ne \hat{x}_k}$ \textbf{then go to 4}}
	
		\State $\mathbf{r_k} = \mathbf{r_{k-1}} - \alpha_{k-1}\mathbf{s_{k-1}}, \quad \mathbf{w'_k} = \mathbf{w_{k-1}} - \alpha_{k-1}\mathbf{u_{k-1}}$
		\State $\nu'_k = \nu_{k-1} - 2\alpha_{k-1} \sigma_{k-1}+\alpha^2_{k-1}\gamma_{k-1}, \quad\beta_k = \nu'_k/\nu_{k-1}$
		\State $\mathbf{p_k} = \mathbf{r_k} + \beta_{k}\mathbf{p_{k-1}}, \quad \mathbf{s_k} = \mathbf{w'_k}+\beta_{k}\mathbf{s_{k-1}}, \quad \mathbf{u_k} = \mathbf{A s_k}, \quad \mathbf{w_k} = \mathbf{A r_k}$
		\State $\mu_k = \langle \mathbf{p_k}, \mathbf{s_k} \rangle,\ \sigma_k = \langle \mathbf{r_k}, \mathbf{s_k} \rangle,\ \gamma_k = \langle \mathbf{s_k}, \mathbf{s_k} \rangle,\ \nu_k = \langle \mathbf{r_k}, \mathbf{r_k} \rangle$
		\State $\alpha_k = \nu_k/\mu_k$
		
		\State $\Delta_{\nu'_k} = |\nu_k - \nu'_k|, \quad {\Delta_{\mathbf{w}'_k} = ||\mathbf{w_{k-1} - w'_{k-1}}||}, \quad\Delta_{\mu'_k} = |\mu_k - \sigma_k|$
		\State $B_{\nu'_k} = \epsilon \, (21+6n)(||\mathbf{r_{k-1}}||^2 + ||\mathbf{r_{k}}||^2)$
		\State $B_{\mathbf{w}'_k} = \epsilon \, ||\mathbf{A}||\, ( (c+3)||\mathbf{r_k}|| + (c+4)||\mathbf{r_{k-1}}|| + (c+2)\,|\alpha_{k-1}|\,||\mathbf{s_{k-1}}|| )$
		\State $B_{\mu'_k} = |\beta_k|\, |\langle  \mathbf{p_{k-1}} , \mathbf{s_k} \rangle| + \epsilon \, ||\mathbf{s_k}||(||\mathbf{r_k}|| + 2|\beta_k|\,||\mathbf{p_{k-1}}|| + n \, (||\mathbf{p_k}|| + ||\mathbf{r_k}||))$
		\NoThenIf{$\Delta_{\nu'_k}>B_{\nu'_k}$ \textbf{or}  {$\Delta_{\mathbf{w}'_{k-1}}>B_{\mathbf{w}'_{k-1}}$} \textbf{or} $\Delta_{\mu'_k}>B_{\mu'_k}$ \textbf{or} $|B_{\mu'_k} - \Delta_{\mu'_k}|/B_{\mu'_k} < T$} 
		\State{\textbf{if} $|B_{\mu'_k} - \Delta_{\mu'_k}|/B_{\mu'_k} < T$ \textbf{then} $T = a\cdot T$}
        \State \textbf{Recover(){, $k=k-3$}}
		\EndIf
		\EndWhile
		\EndProcedure
	\end{algorithmic}
\end{algorithm}
}

We perform one final numerical experiment to investigate the detection reliability of AFT-Pipe-PR-CG and the average number of alarms raised during its runs, the results of which are presented in Table \ref{table3}. The setup of this experiment was very similar to that of the detection experiment in the previous section. The choice of the random problem parameters such as right-hand side or flip iteration was the same. Identical also were the convergence criterion, the initial guess, and that the already calculated variables were utilized for computation of the norms in the detection criteria. The initial value of the threshold was set to $5\mathrm{e}{-1}$. For each of the Pipe-PR-CG variables, excluding $\mathbf{x_k}$, 500 tainted runs were performed, i.e., there were 6500 runs in total for each matrix and choice of $a$. We have decided to include in the experiment only the matrices \textit{1138\_bus} and \textit{nos7}, since for the other matrices from our sample almost no false positives were indicated in the large detection experiment (see Table \ref{table2}).
\par
Runs for the two adaptation parameters $a$ were performed separately, unlike in the case of the experiment in the previous section where the two threshold settings were tested on the same data. The reason for this is that it would be rather difficult to deal with situations when one detection set using the first parameter $a$ does not raise the alarm, but the other detection set using the second parameter $a$ does, and thus it also wants to perform a rollback. Nonetheless, the main purpose of this experiment was not a straight comparison of the two choices of $a$, but rather to investigate whether the introduction of the adaptive threshold refinement causes additional false negative detections, as well as to get a notion of how many alarms there are in average in each run. Additionally, we also obtain information about the number of extra iterations performed due to the recoveries as this is twice the number of alarms.
\par
As  before, first an untainted run to obtain the number of iterations needed to converge for the given right-hand side $\mathbf{b}$ was performed. Then, we executed a tainted run of AFT-Pipe-PR-CG during which it was noted whether the bit flip was detected and corrected (positive cases) or whether the alarm was not raised (special/false negative cases based on the number of iterations to converge) either in the flip iteration or in the following iteration. The positive cases are no longer sorted, since successful detection leads to only two additional iterations, and thus the method always converged within the given limit. The table also contains additional column with the average number of alarms during run.
\par
From Table \ref{table3}, we observe that the introduction of the adaptive threshold refinement did not increase the number of false negative detections, which remains at the same level as it was for the static threshold approach. Moreover, by utilizing the adaptive strategy, we were able to restrict the number of false positive detections to only a handful per run. This is especially impressive for the matrix \textit{nos7}, because of its rather high condition number and the large number of false positive detections observed for it in the detection performance experiment (see Table \ref{table2}). Interesting also may be that the number of average alarms is for the matrix \textit{1138\_bus} very close to one. Nonetheless, this number also includes the negative detections for which there were no alarms at the flip iterations. The AFT-Pipe-PR-CG algorithm is thus strongly reliable and, in combination with a suitable parameters $a$ and $T$, can effectively reduce the number of extra iterations performed due to recovery. 
\begin{table}[h]
\centering
\begin{tabular}{c c c c c c} 
	matrix&$a$&positive& sn & fn & \#alarms \\
	\hline
	\textit{1138\_bus}&$0.5$ & 4192 & 2232 & 4 & 1.057 \\
	&$0.1$&4310 & 2117 & 2 & 1.010 \\
	\hline
	\textit{nos7}&$0.5$&3579 & 2846 & 1 & 13.655\\
	&$0.1$&3487 & 2942 & 1 & 4.784\\
	\hline
\end{tabular}
\renewcommand{\arraystretch}{1}
\caption{Performance of AFT-Pipe-PRCG} 
\label{table3}
\end{table}

\par
As a final remark, we note that the threshold for the relative difference of the $\mu$-gap and the $\mu$-gap bound may also be set based on an estimation of the condition number of $\mathbf{A}$. In the detection experiment we have observed that the higher the condition number of the matrix, the more likely a false positive. As mentioned, the $||\mathbf{A}||$ can be reasonably estimated within few iterations of the Pipe-PR-CG algorithm. Additionally, it is also possible to estimate the condition number \cite{Anorm}, and thus, we could use this information for setting the threshold value. This investigation is left for future research.

\section{Conclusion}\label{sec13}

This article has explored the problem of the detection and correction of silent errors in the Pipe-PR-CG algorithm. 
Our approach is based on the derivation of {finite precision error bounds} for three so-called ``gaps'' between variables which are equal in exact arithmetic. We showed that the violation of these bounds by the computed gaps can be used to detect silent errors in many of the Pipe-PR-CG variables. In order to detect faults in the variables not covered by the three bound violation criteria, a fourth criterion has been constructed, based on monitoring the relative difference between the $\mu$-gap and the $\mu$-gap bound. We then demonstrated that the derived criteria are able to reliably detect the vast majority of silent errors which, if left uncorrected, would significantly impact convergence of the method. In cases when the injected errors remained undetected, the algorithm almost always reaches the stopping criterion without serious delay. 

We then incorporated the derived detection methods along with a recovery procedure into the FT-Pipe-PR-CG algorithm. 
However, it was noted that for some matrices the fault-tolerant algorithm could be significantly slowed down by many extra iterations due to recovery caused by a large number false positive detections. To remedy this, we have proposed the idea of adaptive threshold refinement based on the number of detected alarms during the computation. The resulting adaptive fault-tolerant algorithm, AFT-Pipe-PR-CG, can effectively limit the number of false positives.

{ We note that our approach here requires the computation of additional inner products. We stress that these can be computed in the same single global synchronization point in each iteration, and thus this should not be a significant additional cost in latency-bound regimes. We also note that additional quantities, such as the norm of the input matrix (or norms of preconditioners) appear in the bounds; we again stress that due to the worst case nature of the rounding error analysis, only very rough estimates of these quantities are needed. }

{While we also provided finite precision bounds for the preconditioned case, our experimentation focused on the unpreconditioned Pipe-PR-CG algorithm for simplicity. The derivation of particular bounds that are tight enough to be used effectively to detect silent errors is necessarily highly dependent on the type of preconditioning used and the manner in which the preconditioner is applied. This is a challenge which should be addressed in the future. }

\section*{Acknowledgements}
The first author is supported by GAUK project No. 202722, Charles University Research Centre program No. UNCE/24/SCI/005, and by the European Union (ERC, inEXASCALE, 101075632). Views and opinions expressed are those of the authors only and do not necessarily reflect those of the European Union
or the European Research Council. Neither the European Union nor the granting authority can be held responsible
for them.

\bibliography{main.bbl}%

\begin{thebibliography}{10}

\bibitem{agullo2017hard}
Emmanuel Agullo, Siegfried Cools, Luc Giraud, Alexandre Moreau, Pablo Salas,
  Wim Vanroose, Emrullah~Fatih Yetkin, and Mawussi Zounon.
\newblock Hard faults and soft-errors: possible numerical remedies in linear
  algebra solvers.
\newblock In {\em High Performance Computing for Computational Science--VECPAR
  2016: 12th International Conference, Porto, Portugal, June 28-30, 2016,
  Revised Selected Papers 12}, pages 11--18. Springer, 2017.

\bibitem{CGSoft}
Emmanuel Agullo, Siegfried Cools, Emrullah~Fatih Yetkin, Luc Giraud, Nick
  Schenkels, and Wim Vanroose.
\newblock On soft errors in the conjugate gradient method: Sensitivity and
  robust numerical detection.
\newblock {\em SIAM Journal on Scientific Computing}, 42(6):C336--C358, 2020.

\bibitem{Aupy2017}
Guillaume Aupy, Anne Benoit, Aurlien Cavelan, Massimiliano Fasi, Yves Robert,
  Hongyang Sun, and Bora U{\c{c}}ar.
\newblock {\em Coping with Silent Errors in HPC Applications}, pages 269--292.
\newblock Springer International Publishing, Cham, 2017.

\bibitem{bronevetsky2008soft}
Greg Bronevetsky and Bronis de~Supinski.
\newblock Soft error vulnerability of iterative linear algebra methods.
\newblock In {\em Proceedings of the 22nd Annual International Conference on
  Supercomputing}, pages 155--164, 2008.

\bibitem{Carsonetal}
Erin Carson, Miroslav Rozložník, Zdeněk Strakoš, Petr Tichý, and Miroslav
  Tůma.
\newblock The numerical stability analysis of pipelined conjugate gradient
  methods: Historical context and methodology.
\newblock {\em SIAM Journal on Scientific Computing}, 40(5):A3549--A3580, 2018.

\bibitem{ChenCarson}
Tyler Chen and Erin Carson.
\newblock Predict-and-recompute conjugate gradient variants.
\newblock {\em SIAM Journal on Scientific Computing}, 42(5):A3084--A3108, 2020.

\bibitem{ChGCG}
A.T. Chronopoulos and C.W. Gear.
\newblock s-step iterative methods for symmetric linear systems.
\newblock {\em Journal of Computational and Applied Mathematics},
  25(2):153--168, 1989.

\bibitem{SuiteSparseArticle}
Timothy~A. Davis and Yifan Hu.
\newblock The {U}niversity of {F}lorida sparse matrix collection.
\newblock {\em ACM Trans. Math. Softw.}, 38(1), dec 2011.

\bibitem{ElliottIEEE}
James Elliott, Mark Hoemmen, and Frank Mueller.
\newblock Evaluating the impact of {SDC} on the {GMRES} iterative solver.
\newblock pages 1193--1202, 2014.

\bibitem{elliott2015numerical}
James Elliott, Mark Hoemmen, and Frank Mueller.
\newblock A numerical soft fault model for iterative linear solvers.
\newblock In {\em Proceedings of the 24th International Symposium on
  High-Performance Parallel and Distributed Computing}, pages 271--274, 2015.

\bibitem{GVCG}
P.~Ghysels and W.~Vanroose.
\newblock Hiding global synchronization latency in the preconditioned conjugate
  gradient algorithm.
\newblock {\em Parallel Computing}, 40(7):224--238, 2014.

\bibitem{bitstring}
Scott Griffiths.
\newblock Python module bitstring.
\newblock (Version 4.1).

\bibitem{OGCG}
M.~R. Hestenes and E.~Stiefel.
\newblock Methods of conjugate gradients for solving linear systems.
\newblock {\em Journal of Research of the National Bureau of Standards}, 49(6),
  1952.

\bibitem{martinsson2020randomized}
Per-Gunnar Martinsson and Joel~A Tropp.
\newblock Randomized numerical linear algebra: Foundations and algorithms.
\newblock {\em Acta Numerica}, 29:403--572, 2020.

\bibitem{MeurantCG}
G.~Meurant.
\newblock Detection and correction of silent errors in the conjugate gradient
  algorithm.
\newblock {\em Numerical Algorithms}, 92:869–891, 2023.

\bibitem{Anorm}
G.~Meurant and P.~Tichý.
\newblock Approximating the extreme {R}itz values and upper bounds for the
  {A}-norm of the error in {CG}.
\newblock {\em Numerical Algorithms}, 82:937–968, 2019.

\bibitem{MeurantCray}
Gérard Meurant.
\newblock Multitasking the conjugate gradient method on the {CRAY X-MP/48}.
\newblock {\em Parallel Computing}, 5:267--280, 1987.

\bibitem{SuiteSparse}
The~University of~Florida.
\newblock {SuiteSparse} matrix collection.
\newblock (Last accessed on 2023/11/25).

\end{thebibliography}

\appendix

\section{Statements of the initialization procedures}\label{init}

\begin{algorithm}[h]
	\caption{Initialize (HS-CG)}
	\begin{algorithmic}[1]
		\Procedure {Initialize}{$\mathbf{A}, \mathbf{M}, \mathbf{b}, \mathbf{x_0}$}
		\State $\mathbf{r_0} = \mathbf{b - Ax_0},\quad\nu_0 = \langle \mathbf{\tilde{r}_0}, \mathbf{r_0} \rangle,\quad\mathbf{p_0} = \mathbf{\tilde{r}_0},\quad\mathbf{s_0} = \mathbf{Ap_0},\quad\alpha_0 = \nu_0/\langle \mathbf{p_0}, \mathbf{s_0} \rangle$
		\EndProcedure
	\end{algorithmic}
\end{algorithm}

\begin{algorithm}[h]
	\caption{Initialize (Pipe-PR-CG and variants)}
	\begin{algorithmic}[1]
		\Procedure {Initialize}{$\mathbf{A}, \mathbf{M}, \mathbf{b}, \mathbf{x_0}$}
		\State $\mathbf{r_0} = \mathbf{b - Ax_0},\quad  \mathbf{\tilde{r}_0} = \mathbf{M^{-1}r_0},\quad\mathbf{p_0} = \mathbf{\tilde{r}_0},\quad\mathbf{s_0} = \mathbf{Ap_0},\  \mathbf{\tilde{s}_0} = \mathbf{M^{-1}s_0}$
		\State $\mathbf{u_0} = \mathbf{A\tilde{s}_0} ,\quad  \mathbf{\tilde{u}_0} = \mathbf{M^{-1}u_0},\quad\mathbf{w_0} = \mathbf{A\tilde{r}_0} ,\quad  \mathbf{\tilde{w}_0} = \mathbf{M^{-1}w_0} $
		\State $\sigma_0 = \langle \mathbf{\tilde{r}_0}, \mathbf{s_0} \rangle,\quad\gamma_0 = \langle \mathbf{\tilde{s}_0}, \mathbf{s_0} \rangle,\quad\nu_0 = \langle \mathbf{\tilde{r}_0}, \mathbf{r_0} \rangle,\quad\alpha_0 = \nu_0/\langle \mathbf{p_0}, \mathbf{s_0} \rangle$
		\EndProcedure
	\end{algorithmic}
\end{algorithm}

\section{Detection Criteria Experiments}\label{AppendixFigures}
Here we present additional figures displaying the behavior of quantities utilized in the presented detection criteria. We use matrix ($\mathbf{A}$) \textit{nos7}, right-hand side $\mathbf{b}=\mathbf{Ae}$, and $x_0$ as a vector of all zeros. Flips occurred in the 15th bit (for vector variables in the position 100).
\begin{figure}[h]
	\makebox[\textwidth][c]{\includegraphics[clip, trim=0cm 0cm 0cm 1cm, width=.7\textwidth]{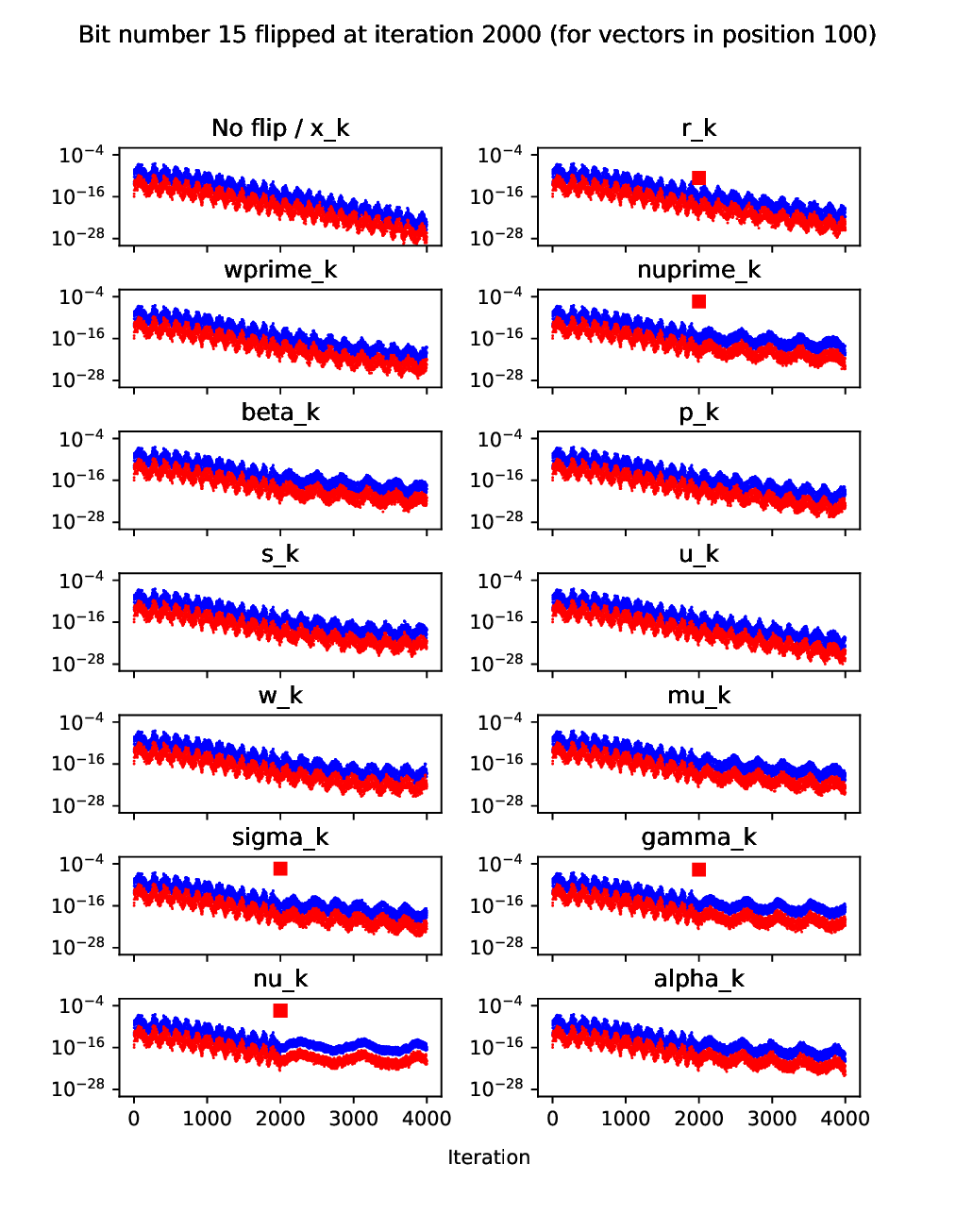}}
	\caption{$\nu$-gap (red) and $\nu$-gap bound (blue) graph, matrix \textit{nos7}}
	\label{FigB1}
\end{figure}
\begin{figure}[h]
	\makebox[\textwidth][c]{\includegraphics[clip, trim=0cm 0cm 0cm 1cm, width=.7\textwidth]{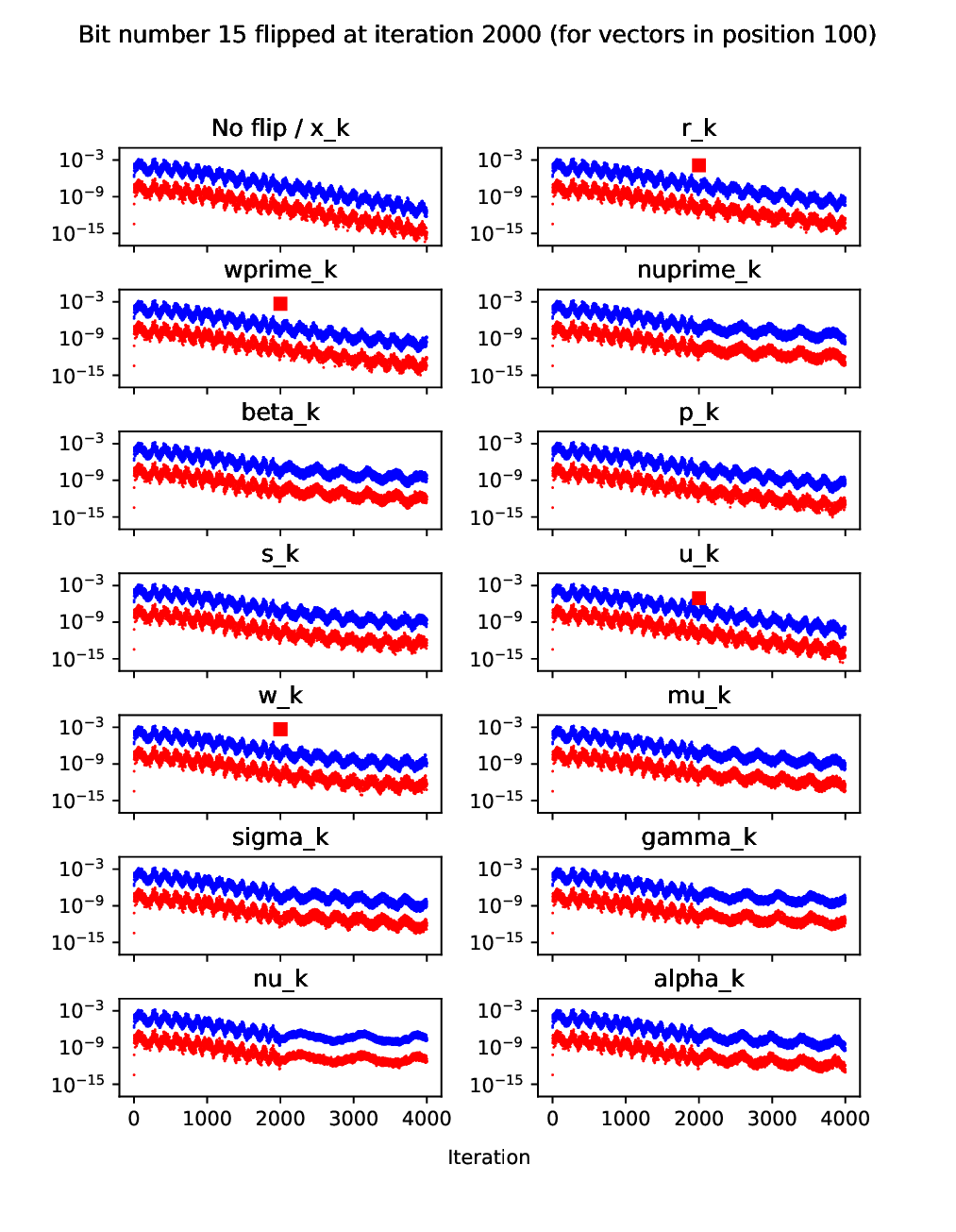}}
	\caption{$\mathbf{w}$-gap (red) and $\mathbf{w}$-gap bound (blue) graph, matrix \textit{nos7}}
	\label{FigB2}
\end{figure}
\begin{figure}[h]
	\makebox[\textwidth][c]{\includegraphics[clip, trim=0cm 0cm 0cm 1cm, width=.7\textwidth]{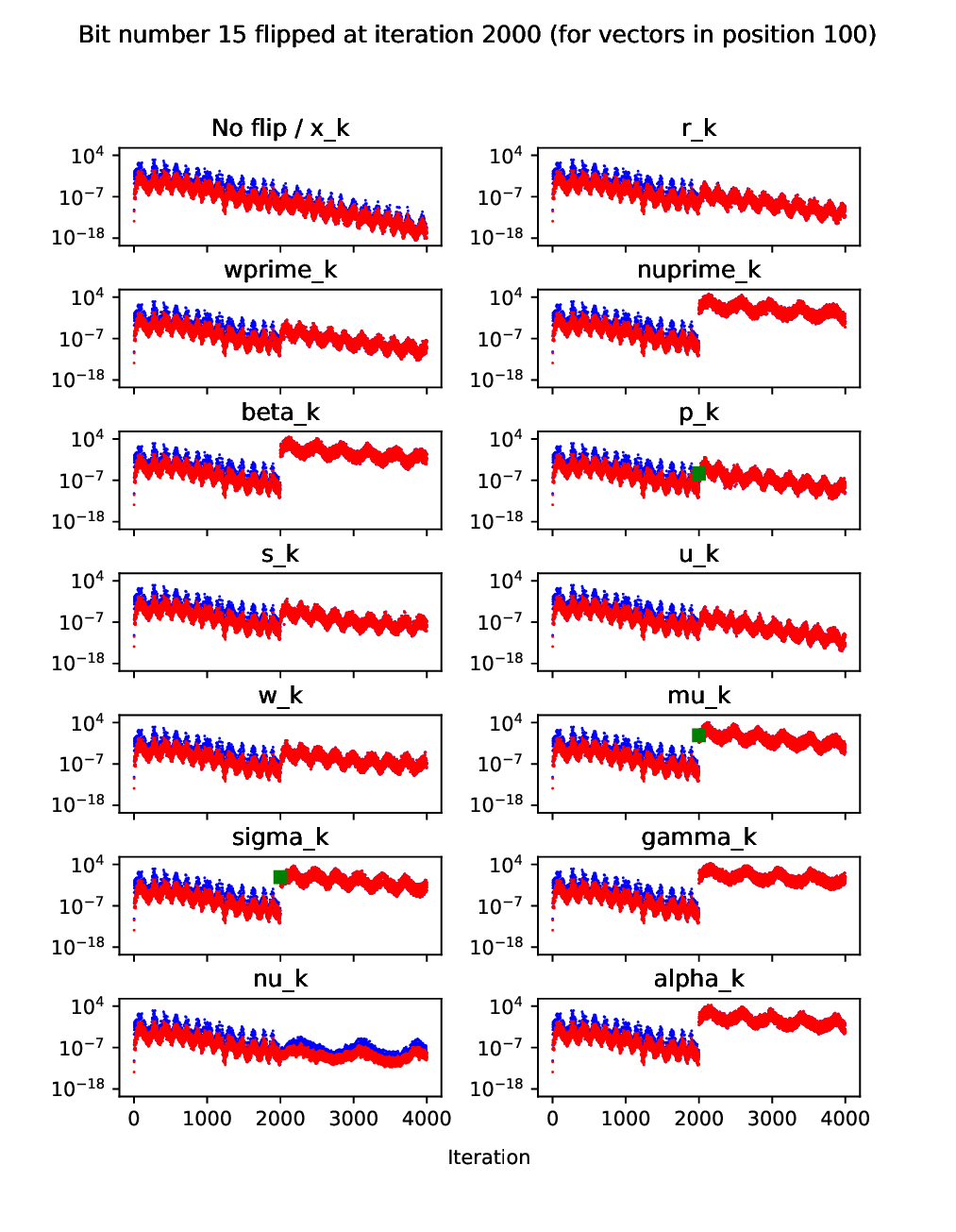}}
	\caption{$\mu$-gap (red) and $\mu$-gap bound (blue) graph, matrix \textit{nos7}}
	\label{FigB3}
\end{figure}

\begin{figure}[h]
	\makebox[\textwidth][c]{\includegraphics[clip, trim=0cm 0cm 0cm 1cm, width=.7\textwidth]{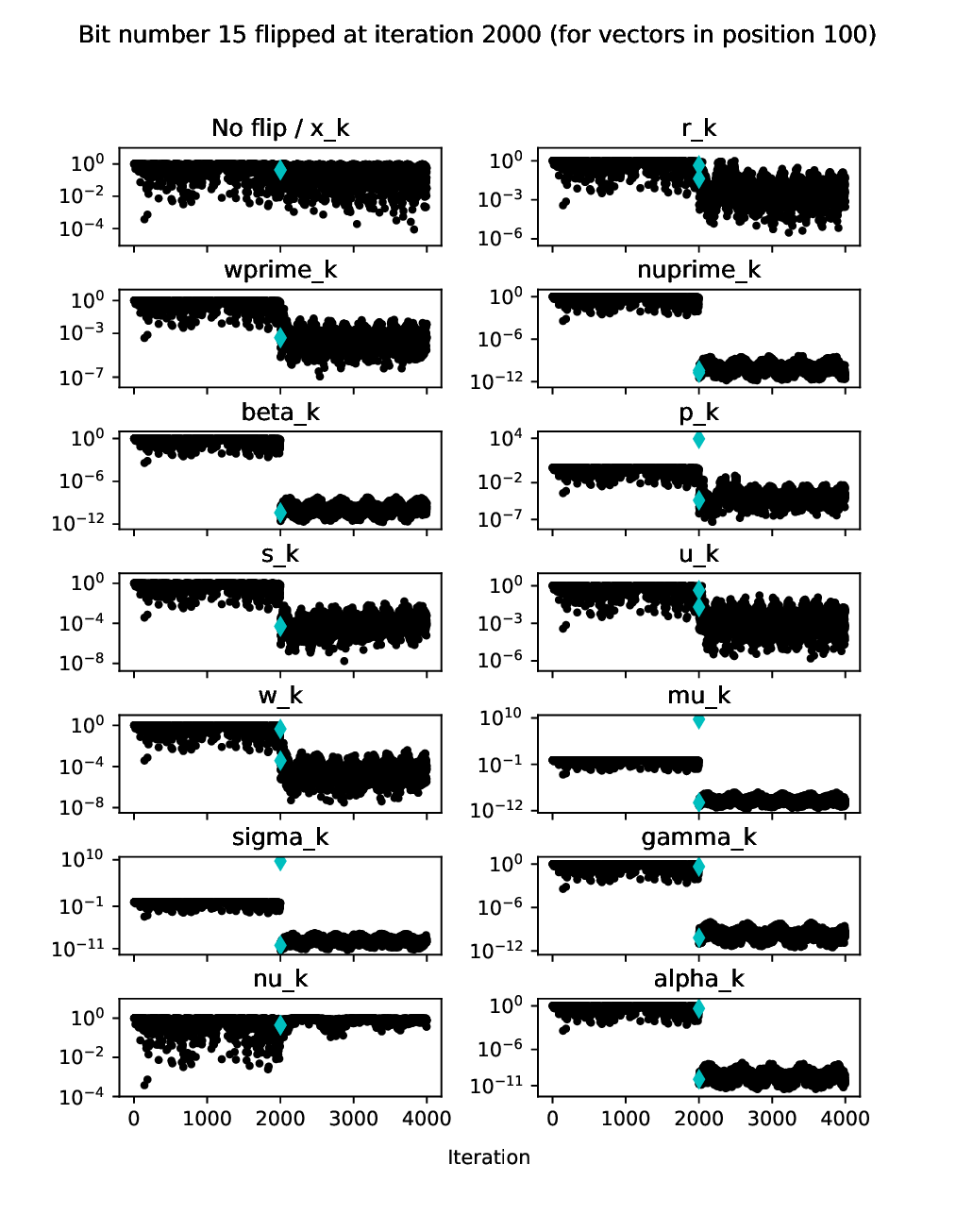}}
	\caption{Relative $\mu$-gap/bound difference, $|B_{\mu'_k} - \Delta_{\mu'_k}|/B_{\mu'_k}$, matrix \textit{nos7}}
     \label{FigB4}
\end{figure}

\end{document}